\documentstyle[12pt,amssymb,pb-diagram]{amsart}
\outer\def\proof{\smallbreak\noindent{\bf Proof.}\enspace}

\def\endproof{\qed}

\def\mmu{{\protect\boldsymbol{\mu}}}

\def\sm{_{\rm sm}}

\def\sh{^{\rm sh}}

\newcommand{\mynewtheorem}[3]{\newenvironment{#1}{\subsection{#2} #3}{}}

\mynewtheorem{th}{Theorem}{\it}
\mynewtheorem{claim}{Claim}{\it}
\mynewtheorem{plem}{Purity Lemma}{\it}	
\mynewtheorem{lem}{Lemma}{\it}		\mynewtheorem{situ}{Situation}{\it}
\mynewtheorem{cor}{Corollary}{\it}	\mynewtheorem{prp}{Proposition}{\it}
\mynewtheorem{dfn}{Definition}{\rm}	\mynewtheorem{notation}{Notation}{\rm}
\mynewtheorem{example}{Example}{\rm}	\mynewtheorem{conj}{Conjecture}{\rm}
\mynewtheorem{prob}{Problem}{\rm}	\mynewtheorem{question}{Question}{\rm}
\mynewtheorem{discuss}{}{\rm}
\mynewtheorem{rem}{Remark}{\rm}

 	\newcommand{\dar}{\downarrow}
\newcommand{\down}{\downarrow}	\newcommand{\into}{\hookrightarrow}
	
\newcommand{\bfp}{{\bold{P}}} 	\newcommand{\bfc}{{\bold{C}}}
\newcommand{\bfq}{{\bold{Q}}}	
\newcommand{\bfz}{{\bold{Z}}}	
	\newcommand{\bfa}{{\bold{A}}}
\newcommand{\co}{{\cal{O}}}

\def\mm{{\frak m}}

\def\mapright#1{\mathrel
{\smash{\mathop{\longrightarrow} \limits^{#1}}}}

\def\double{{\raise-1pt\hbox{$\scriptstyle\to$} \atop
\raise3pt\hbox{$\scriptstyle\to$}}}

\newcommand{\spec}{{\operatorname{Spec\ }}}
\newcommand{\sing}{{\operatorname{Sing }}}
\newcommand{\spf}{{\operatorname{Spf\ }}}

\newcommand{\Aut}{{\operatorname{Aut}}}
\newcommand{\aut}{{\operatorname{Aut}}}

\newcommand{\isom}{\operatorname{Isom }}

\newcommand{\CModuli}{{\overline{\bold{M}}}}
\newcommand{\m}{{\overline{\bold M}_{\gamma, \nu}}}
\newcommand{\Moduli}{\overline{\cal M}}
\newcommand{\M}{{\overline{\cal M}_{\gamma, \nu}}} 
\newcommand{\F}{{\cal F}}
\newcommand{\G}{{\cal G}}
\newcommand{\C}{{\cal C}}
\newcommand{\X}{{\cal X}}
\newcommand{\Y}{{\cal Y}}
\newcommand{\g}{{\frak g}}

\newcommand{\curveuparrow}{%
\begin{picture}(4,10)%
\put(0,10){\vector(-1,2){0}}%
\qbezier(0,0)(4,5)(0,10)
\end{picture}}

\newcommand{\rest}{|_}

\renewcommand{\setminus}{%
\begin{picture}(6,10)\qbezier(1,5.5)(3,4.5)(5,3.5)\end{picture}}

\def\into{\hookrightarrow}

\def\cursrhom{\mathop{{\rm R}\mit{\cal H}om}\nolimits}

\def\tototi{\mathbin{\mathop{\otimes}\limits^{\raise-1pt\hbox
{$\scriptscriptstyle {\rm L}$}}}}

\def\ext{\mathop{\rm Ext}\nolimits}

\def\PP{{\cal F}^d_g(\gamma,\nu)}

\def\pp{{\bold F}^d_g(\gamma,\nu)}

\def\K{{\cal K}^d_g\bigl(\,\m\bigr)}

\def\k{{\bold K}^d_g\bigl(\,\m\bigr)}

\begin{document}
\begin{center}
{\bf COMPLETE MODULI FOR FIBERED SURFACES} \end{center}
{\bf Preliminary version}\\[2mm]
{Dan Abramovich\footnote{Partially supported by NSF grant DMS-9503276 and
an Alfred P. Sloan Research Fellowship.} \\ 
Department of
Mathematics, Boston University\\ 111 Cummington Street, Boston, MA 02215,
USA \\ {\tt abrmovic@@math.bu.edu}}\\[2mm]
{Angelo Vistoli\\
Dipartimento di Matematica,
Universit\`a di Bologna \\
Piazza di Porta San Donato 5,
40127 Bologna, Italy \\
{\tt vistoli@@dm.unibo.it}}\\[2mm]
\today


\section{Introduction}\label{intro}

All the schemes with which we work will be schemes over the field ${\bfq}$
of rational numbers.\footnote{The theory below could be extended to the
case that the characteristic is large with respect to the degrees and
genera of the curves involved, but we will not pursue this here.}

\subsection{The problem}
Fix a base field $k$ of characteristic 0, and let $C$ be a smooth,
projective, geometrically integral curve of genus $g$. By a {\em fibered
surface over $C$} we mean a morphism $X\to C$, with sections
$\sigma_1,\ldots,\sigma_\nu:C\to X$ forming a family of stable
$\nu$-pointed curves of some genus $\gamma$. (This notion will be
generalized below for singular $C$.)

A fibered surface naturally corresponds to a morphism $C\to \M$, into the
moduli stack of stable $\nu$-pointed curves of genus $\gamma$. In
\cite{vistoli1}, the second author showed that $Hom(C,\M)$ has the
structure of a Deligne - Mumford stack. This stack is clearly not complete
in general. A natural question to ask is, can one complete it in a
meaningful way? Moreover, can one do this as the source curve $C$ itself
moves in a family?

\subsection{Our approach} One construction which may spring to mind is that of
{\em stable maps}. 
Recall that if $Y\subset \bfp^r$ is a projective variety, there is a complete
Deligne - 
Mumford stack $\Moduli_g(Y,d)$ parametrizing Kontsevich stable maps of
degree $d$, from curves of genus $g$ to $Y$ (see  \cite {k}, \cite{bm},
\cite{fp}). In order to reduce the possibility of later confusion, we will use
the notation ${\cal K}_g^d(Y)$ rater than $\Moduli_g(Y,d)$ for this stack.
It admits a projective coarse moduli space ${\bold K}_g^d(Y)$.

For instance, if $\gamma=0$, the moduli stack $\Moduli_{0,\nu}$ is actually
a projective variety. Fixing a projective embedding, we have natural
complete Deligne - Mumford stacks ${\cal K}_g^d(\Moduli_{0,\nu})$
parametrizing families of stable $\nu$-pointed curves of genus 0 over nodal
curves of genus $g$, with a suitable stability condition. We can think of
$d$, the degree of the image of $C$ via the fixed projective embedding, as
an additional measure of complexity of the family of curves. One is led to
ask, is there a complete Deligne - Mumford stack of stable maps ${\cal
K}_g^d(\Moduli_{\gamma,\nu})$ in general?

We need to define our terms: given a nodal curve $C$, define a morphism $f:C\to
\M$ to be a stable map of degree $d$, if once projected into the coarse moduli
space, we obtain a stable map $f':C\to \m$ of degree $d$: 

$$\begin{diagram}
\node{C}\arrow{e,t}{f}\arrow{se,l}{f'}\node{\M}\arrow{s}\\
\node[2]{\m}
\end{diagram}$$

As we will see later, these stable maps are parametrized by a Deligne -
Mumford stack. However, this stack fails to perform its first goal: it is
not always complete! This can be seen via the following example:

Consider the case $\gamma=1, \nu=1$ of elliptic curves. It is well known
(and easily follows from the formula of Grothendieck - Riemann - Roch) that
given a family of stable elliptic curves $E\to C$ over a curve $C$, the
degree of the corresponding moduli map $j_E:C\to \CModuli_{1,1}\simeq
\bfp^1$ is divisible by 12. However, it is very easy to construct a family
of moduli maps $j_t:C_t\to \bfp^1,\quad t\neq 0$,
where when $t\to 0$ the curve $C_t$ breaks in two components, and the
degrees of the limit map $j_0$ on these components are not divisible by 12.
Thus $j_0$ is not a moduli map of a family of stable curves.

\subsection{Compactifying the space of stable maps} Our main goal in this paper
is to correct this deficiency. In order to do 
so, we will enlarge the category of stable maps into $\M$. The source curve
$\C$ of a new stable map $\C \to \M$ will acquire an orbispace structure at
its nodes. Specifically, we endow it  with the structure of a Deligne-Mumford
stack.

To see how these come about, consider again the example of a one-parameter
family of elliptic curves sketched above. Let $C\to S$ be the one-parameter
family of base curves, and let $C\sm$ be the smooth locus of this morphism.
A fundamental purity lemma (see \ref{purity-lemma} below)  will show,
that after a suitable base change we can extend the family $E$ of stable
elliptic curves 
over $C\sm$. On the other hand, if $p\in C$ is a node, then on an \'etale
neighborhood $U$ of $p$, the
curve $C$ looks like $$uv = t^r,$$
where $t$ is the parameter on the base. By taking $r$-th roots:
$$u = u_1^r;\, v = v_1^r$$ we have a nonsingular covering $V_0\to U$ where
$u_1v_1 = t$. The fundamental purity lemma applies to $V_0$, so the
pullback of $E$ to $V_0$ extends over all of $V_0$. There is a minimal
intermediate covering $V_0\to V\to U$ such that the family extends already
over $V$. This $V$ gives the orbispace structure $\C$ over $C$.

The precise definitions will follow in section \ref{sec:fibered-surfaces}.
We can define stable maps $f:\C \to \M$ where $\C$ is a Deligne-Mumford stack
and $f$ is a representable morphism satisfying certain properties; such objects
naturally form a 2-category. Our main theorem states that this 2-category is
equivalent to a 1-category $\PP$, which we call {\em the category of fibered
surfaces};  this category forms a complete Deligne-Mumford stack
(over schemes over $\bfq$), admitting a projective coarse moduli space. Thus
the original moduli problem has a solution of the same nature as that of stable
curves.

We provide an explicit description of the category $\PP$ of fibered surfaces in
terms of charts and atlases over schemes, in analogy to Mumford's treatment of
$\bfq$-varieties in \cite{MumTowards}. Thus given a stable map $\C \to \M$ we
have a stable pointed curve $\X \to \C$, which can be described via an atlas of
charts over the associated coarse moduli scheme $X \to C$.

\subsection{Comparison with Alexeev's work}
The latter, coarse object $X \to C$, which we call a {\em coarse fibered
surface}, has another interpretaion: the associated morphism $X \to \m$ turns
out to be a stable map in the sense of Alexeev (see \cite{alexMaps}). Alexeev
has shown the existence of complete moduli of (smoothable) surface stable maps,
using the theory of semi-log-canonical surfaces developed in \cite{ksb} and
\cite{alex}. In the last section we show that in the case of fibered surfaces,
Alexeev's approach coincides with the space of (smoothable) coarse fibered
surfaces. We sketch a proof of the existence of the latter space (even as a
stack)  which is
independent of Alexeev's boundedness result.

\subsection{Natural generalizations in forthcoming work} It should be evident
from our work that 
this approach should apply to stable maps $\C \to \Moduli$ into any
Deligne-Mumford stack which admits a projective coarse moduli scheme; moreover,
the source ``curve'' $\C$ should be allowed to be ``pointed'' as well. We are
now in the process of writing up the necessary details, circumventing the
various shortcuts we have taken in this paper. Various applications
promise to be of interest; e.g. the case $\Moduli = B{\frak S}_n$, the
classifying space of the symmetric group on $n$-letters, provides us with a
natural desingularization of the moduli stack of admissible covers (whose
coarse moduli scheme was introduced
by Harris and Mumford \cite{hm}; see  \cite{moch} for a rigorous construction
using logarithmic structures). The case $\Moduli =B((\bfz/n\bfz)^{2g})$ can be
used to give a moduli interpretation for Mumford's compactification of the
moduli space of curves with level structures (see \cite{MumTowards}), and the
refinements by Looijenga \cite{l}
and Pikaart - de Jong \cite{p-dj}. Also, the construction of stable maps is
amenable to an inductive approach, allowing us to construct a moduli space of
``plurifibered varieties'', giving  complete moduli in arbitrary
dimension for varieties with ``plurifibrations''. This is closely related to
work of Mochizuki \cite{mochfib1,mochfib2}, where  the existence of
 level structure is assumed.

\subsection{Gromov - Witten invariants} Originally, the Kontsevich
spaces of stable maps were introduced for the purpose of defining
Gromov - Witten invariants.  It seems likely that, using our
construction,  one could extend the formalism of
[BF] and [B] of virtual fundamental classes to the
case of stable maps into an arbitrary Deligne Mumford stack $\Moduli$ admitting
a projective  coarse moduli space. This should allow one to define
Gromov - Witten invariants and quantum cohomology in this generality, which may
have interesting application for specific choices of stacks $\Moduli$.

\subsection{Acknowledgements} We would like to thank Johan de
Jong and Rahul Pandharipande for crucial discussions of ideas in this paper.

\section{The purity lemma}

There are several results in the litterature, which give conditions under which
a family of curves can fail to be stable only in pure codimension 1 (see
\cite{mb}, \cite{dj-o}). For our purposes, the following case will be most
useful: 
\begin{plem}\label{purity-lemma} Let $\Moduli$ be a
separated 
Deligne - Mumford 
stack, $\Moduli\to
\CModuli$ the coarse moduli space. Let $X$ be a separated scheme of
dimension 2 satisfying Serre's condition $S_2$. Let $P\subset X$ be a finite
subset consisting of closed points, $U=X\setminus P$. Assume that the local
fundamental groups of $U$
around the points of $P$ are trivial.

Let $f:X \to \CModuli$ be a morphism. Suppose there is a lifting
$\tilde{f}_U:U \to \Moduli$:
\begin{equation}\begin{diagram} \node[3]{\Moduli}\arrow{s} \\
\node{U}\arrow{ene,t}{\tilde{f}_U}\arrow{e} 
\node{X}\arrow{e,t}{f} \node{\CModuli} \end{diagram}
\end{equation}

Then the lifting extends uniquely to $X$:
\begin{equation}\begin{diagram} \node[3]{\Moduli}\arrow{s} \\
\node{U}\arrow{ene,t}{\tilde{f}_U}\arrow{e} 
\node{X}\arrow{e,t}{f}\arrow{ne,t,1,..}{\tilde{f}} \node{\CModuli}
\end{diagram}
\end{equation}
\end{plem}

\proof
By descent theory the problem is local in the \'etale topology, so we may
replace $X$ and $\CModuli$ with the spectra of their strict henselizations
at a geometric point; then we can also assume that we have a universal
deformation space
$V\to \Moduli$ which is {\em finite}. Now $U$ is the complement of the
closed point, $U$ maps to
$\Moduli$, and the pullback of $V$ to $U$ is finite and \'etale, so it has
a section, because $U$ is simply connected; consider the corresponding map
$U\to V$. Let $Y$ be the scheme-theoretic closure of the graph of this map
in $X\times_\CModuli V$. Then $Y\to X$ is finite and is an isomorphism on $U$.
Since $X$ satisfies  $S_2$, the morphism $Y\to X$ is an isomorphism. \qed

\begin{cor}
Let $X$ be a smooth surface, $p\in X$ a closed point with complement $U$.
Let $X\to \CModuli$ and $\, U\to \Moduli$ be as in the purity lemma. Then
there is a unique lifting $X \to \Moduli$.
\end{cor}

\begin{cor}
Let $X$ be a normal crossings surface, namely a surface which is
\'etale locally isomorphic to $\spec k[u,v,t]/(uv)$. Let $p\in X$ a closed
point with complement $U$. Let $X\to \CModuli$ and $\, U\to \Moduli$ be as
in 
the purity lemma. Then there is a unique lifting $X \to \Moduli$.
\end{cor}

\proof
In both cases $X$ satisfies condition $S_2$ and the local fundamental group
around $p$ is trivial, hence the purity lemma applies. \endproof

\section{Group actions on nodal curves}

Fix two nonnegative integers $\gamma$ and $\nu$. An important ingredient in the
theory will be formed of data as follows:
\begin{enumerate}
\item a diagram $$\begin{array}{c}Y \\ \down\curveuparrow \\ V \\ \down \\ S
 \end{array}$$ 
 (we think of $Y$ and $V$ as $S$-schemes);
\item the morphism $\rho\colon Y\to V$ comes with sections $\tau_1\colon V\to
Y,\ldots, \tau_\nu\colon V\to Y$, forming a family of
stable $\nu$-pointed curves of genus $\gamma$ on $V$.
\end{enumerate} 
Given a finite group $\Gamma$,  by {\em an action of $\Gamma$ on $\rho$} we
mean a pair of actions 
of $\Gamma$ on $V$ and $Y$ as $S$-schemes such that the morphisms $\rho$,
$\tau_1,
\ldots,\tau_\nu$ are $\Gamma$-equivariant. Such an action induces a
morphism $Y/\Gamma \to V/\Gamma$ together with sections $V/\Gamma\to
Y/\Gamma$.

We adopt the convention that all families of curves are assumed to be {\it of
finite presentation} over the base.

\begin{dfn} Let $\rho\colon Y\to V$ be a family of $\nu$-pointed curves of
genus $\gamma$, $\Gamma$ a finite group. An action of $\Gamma$ on $\rho$ is
called\/ {\em essential} if no nontrivial element of $\Gamma$ leaves a
geometric fiber of $\rho$ fixed.
\end{dfn}

Another way to state this condition is to require that if $\gamma$ is an
element of $\Gamma$ leaving a geometric point $v_0$ of $V$ fixed, then
$\gamma$ acts nontrivially on the fiber of $\rho$ on $v_0$.

The following lemma will help us replace an action of  $\Gamma$ on $\rho$  by
an essential action on a related family of curves.

\begin{lem}\label{QuotientIsNodal} Let $V$ be the spectrum of a
 local ring $R$,
$Y\to V$ a 
flat family of nodal curves, $\Gamma_0$ a group acting compatibly on $Y$ and
$V$. Suppose that the action of $\Gamma_0$ on the residue field $k$ of $R$ and
on the fiber
$Y_0$ of $Y$ over $\spec k$ is trivial. Then $Y/\Gamma_0\to V/\Gamma_0$ is
again a flat family of nodal curves.
\end{lem}

\begin{rem} It is easy to give an example showing that this fails in
positive characteristic.
\end{rem}

\proof We need to show that the map $\rho\colon Y/\Gamma_0\to V/\Gamma_0$ is
flat, and the natural map $Y_0\to \rho^{-1}(\spec k)$ is an isomorphism.
Replacing $R$ with the completion of its strict henselization we may assume
that $R$ is complete and $k$ is algebraically closed. Choose a rational point
$p\in X(k)$, and let $M$ be the completion of the local ring of $Y$ at
$p$. Write 
$M\otimes_R k$ as a quotient of a power series algebra $k[[t_1,\ldots,
t_r]] = k[[{\bf t}]]$, and lift the images of the $t_i$ in $M\otimes_R k$ to
invariant elements in $M$. We get a surjective homomorphism $R[[{\bf
t}]]\to M$ which is equivarinant, if we let $\Gamma_0$ act on $R[[{\bf
t}]]$ leaving the $t_i$ fixed. Then the result follows from the lemma below.

\begin{lem} Let $M$ be a finitely generated $R[[{\bf
t}]]$-module which is flat over $R$, such that $\Gamma_0$ acts
trivially on $M\otimes_R k$. Then $M^{\Gamma_0}$ is flat
over $R^{\Gamma_0}$, and the natural homomorphism
$M^{\Gamma_0}\otimes_{R^{\Gamma_0}}R\to M$ is an isomorphism.
\end{lem}

\proof We have $R[[{\bf t}]]^{\Gamma_0} = R^{\Gamma_0}[[{\bf
t}]]$, and from this
we see that the statement holds when $M$ is free over $R[[{\bf
t}]]$. In general, take a finite set of generators of the finite
$k[[{\bf t}]]$-module
$M\otimes_R k$ and lift them to a set of invariant
generators of $M$. We obtain an equivariant surjective
homomorphism $F\to M$, where $F$ is a free $R[[t]]$-module. Let  $K$ be the  
kernel: $$0 \to K \to F \to M \to 0.$$ 
Note that by the finite presentation assumption on $Y \to V$, $K$ is finitely
generated. Since $M$ is flat we get $Tor_i(K,\cdot) 
= Tor_{i+1}(M,\cdot) = 0$, thus $K$ is also flat. Moreover, from the exact 
sequence  $$0 \to K\otimes_Rk \to F \otimes_Rk\to M\otimes_R k\to 0$$ we get
that  $\Gamma_0$ acts trivially on 
$K\otimes_Rk$. We obtain a commutative diagram with exact
rows
$$
\begin{array}{ccccccccc}
&&K^{\Gamma_0}\otimes_{R^{\Gamma_0}}R& \longrightarrow &
F^{\Gamma_0}\otimes_{R^{\Gamma_0}}R& \longrightarrow &
M^{\Gamma_0}\otimes_{R^{\Gamma_0}}R& \longrightarrow & 0\\
&&\down{}&&\down{}&&\down{}&&\\
0& \longrightarrow &K& \longrightarrow &F&\longrightarrow
&M&\longrightarrow &0
\end{array}
$$
where the middle column is an isomorphism. It follows that the
right column is surjective; since $K$ is also flat and
$\Gamma_0$ acts trivially on $K\otimes_R k$ we can apply the same argument to
$K$, and it follows that
the left column is also surjective. Therefore  the right column is
also injective, as desired.

Applying the same argument to the module $K$, we have that the left column is
also an isomorphism, and the first arrow in the top row is injective; the
flatness of 
$M^{\Gamma_0}$ over $R^{\Gamma_0}$ then follows from
Grothendieck's local criterion of flatness.\endproof

The following lemma will help us identify a situation where an essential action
is free: 

\begin{lem}\label{QuotientDropsGenus} Let $G$ be a group of automorphisms of a
stable $\nu$-pointed 
curve $X$ of genus $\gamma$ over an algebraically closed field, and assume that
the quotient $X/G$ is also a stable $\nu$-pointed curve of genus $\gamma$. Then
$G$ is trivial.
\end{lem}

\proof If $C$ is a complete curve we denote by ${\rm CH}^i(C)$ the
Chow group of classes of cycles of codimension $i$ on $C$; so ${\rm
CH}^0(C)$ is a direct sum of $n$ copies of {\bf Z}, where $n$ is the number
of irreducible components of $C$, while ${\rm CH}^1(C)$ is the divisor
class group of $C$. If $\cal F$ is a coherent sheaf on $C$ we denote by
$\tau^C({\cal F})\in {\rm CH}(C)$ its Riemann--Roch class, and by
$\tau^C_i({\cal F})\in {\rm CH}^i(C)$ its component of codimension $i$;
$\tau^C_0({\cal F})$ is the cycle associated with ${\cal F}$.

Set $Y = X/G$, and call $\pi\colon X\to Y$ the projection. Call $\Delta_X$
and $\Delta_Y$ the divisors corresponding to the distinguished points in
$X$ and $Y$ respectively, $\omega_X$ and $\omega_Y$ the dualizing sheaves.

A local calculation reveals that $\chi\bigl(\omega_X(\Delta_X)\bigr) \ge 
\chi\bigl(\pi^*\omega_Y(\Delta_Y)\bigr)$. Indeed, for this calculation we may
replace $X $ and $Y$ by their normalization, where we add markings on the
normalized curves at the points above the nodes. Then 
the natural pullback $\pi^*\omega_{Y^{nor}}\to
\omega_{X^{nor}}$  extends to a
homomorphism of invertible sheaves $\pi^*\omega_{Y^{nor}}(\Delta_{Y^{nor}})\to
\omega_{X^{nor}}(\Delta_{X^{nor}})$ which is injective, and whose cokernel is
supported on  
the points of ramification of $\pi$ not in $\Delta_{X^{nor}}$.

 Thus we get an
inequality of Euler characteristics of sheaves:

\begin{eqnarray*} \chi\bigl(\omega_X(\Delta_X)\bigr)&\ge&
\chi\bigl(\pi^*\omega_Y(\Delta_Y)\bigr)\\ &=&
\chi\bigl(\pi_*\pi^*\omega_Y(\Delta_Y)\bigr)\\ & =&
\chi\bigl(\omega_Y(\Delta_Y)\otimes\pi_*{\cal O}_X\bigr) \\ &=&
\chi\bigl(\omega_Y(\Delta_Y)\bigr) +
\chi\bigl(\omega_Y(\Delta_Y)\otimes(\pi_*{\cal O}_X/{\cal O}_Y)\bigr).
\end{eqnarray*}

By hypothesis
$$
\chi\bigl(\omega_X(\Delta_X)\bigr) =
\gamma + \nu - 1 =
\chi\bigl(\omega_Y(\Delta_Y)\bigr).
$$
{} From this, together with the theorem of Grothendieck--Riemann--Roch, we
obtain 

\begin{eqnarray*} 0&\ge&
\chi\bigl(\omega_Y(\Delta_Y)
\otimes (\pi_*{\cal O}_X/{\cal O}_Y)\bigr) \\& = &\int_Y {\rm
ch}\bigl(\omega_Y(\Delta_Y)\bigr)
\tau^Y(\pi_*{\cal O}_X/{\cal O}_Y)\\
&=& \int_Y {\rm
c}_1\bigl(\omega_Y(\Delta_Y)\bigr)\cap\tau^Y_0(\pi_*{\cal O}_X/{\cal O}_Y) +
\int_Y \tau^Y_1(\pi_*{\cal O}_X/{\cal O}_Y)\\ &=& \int_Y {\rm
c}_1\bigl(\omega_Y(\Delta_Y)\bigr)\cap(\pi_*[X]-[Y]) + \int_X
\tau^X_1({\cal O}_X) - \int_Y\tau^Y({\cal O}_Y)\\ &=& \int_Y {\rm
c}_1\bigl(\omega_Y(\Delta_Y)\bigr)\cap(\pi_*[X]-[Y]) + \chi({\cal O}_X) -
\chi({\cal O}_Y)\\
&=& \int_Y {\rm
c}_1\bigl(\omega_Y(\Delta_Y)\bigr)\cap(\pi_*[X]-[Y]). \end{eqnarray*}

But $\omega_Y(\Delta_Y)$ is an ample line bundle on $Y$ and $\pi_*[X]-[Y]$
is an effective cycle; the only possibility is that $\pi_*[X]-[Y] = 0$.
This implies that $\pi$ is birational, so $G$ is trivial.\endproof

\section{Fibered surfaces: definitions}\label{sec:fibered-surfaces}

If $C\to S$ is a flat family of nodal curves, we denote by $C\sm$ the open
subscheme of $C$ consisting of points where the morphism $C\to S$ is smooth.

\begin{dfn} Let $C\to S$ be a flat (not necessarily proper) family of nodal
curves,
$X\to C$ a proper morphism with one dimensional fibers, and
$\sigma_1,\ldots,\sigma_\nu\colon C\to X$
sections of $\rho$. We will say that
$$\begin{array}{c}X \\ \down \\ C \\ \down \\ S \end{array}$$ is a {\em
family of generically fibered surfaces} if $X$ is flat over $S$, and the
restriction of $\rho$ to $C\sm$ is a flat family of stable pointed curves.
If $S$ is the spectrum of a field we will refer to $X\to C$ as a
generically fibered surface. \end{dfn}

\begin{rem} Notice that we do not require the morphism $X\to C$ to be flat.
\end{rem}

\begin{dfn}
A {\em chart} $(U,Y\to V,\Gamma)$
for a family of generically fibered surfaces $X\to C\to S$ consists of a
diagram 
$$\begin{diagram}
\node{Y}\arrow{s}\arrow{e}
	\node{X\times_CU}\arrow{s}\arrow{e} 
		\node{X}\arrow{s} \\
\node{V}\arrow{e}\arrow{ese}
	\node{U}\arrow{e}\arrow{se} 
		\node{C}\arrow{s} \\
		\node[3]{S}
\end{diagram}
$$
together with a group action $\Gamma\subset \aut_S(Y\to V)$ satisfying:
\begin{enumerate}
\item The morphism $U\to C$ is \'etale;
\item $V\to S$ is a flat (but not necessarily proper) family of nodal curves;
\item $\rho\colon Y\to V$ is a flat family of stable
$\nu$-pointed curves of genus $\gamma$,
\item the action of $\Gamma$ on $\rho$ is essential; \item we have
isomorphisms of
$S$-schemes $V/\Gamma\simeq U$ and $Y/\Gamma\simeq U\times_C X$ compatible
with the projections $Y/\Gamma \to V/\Gamma$ and $U\times_C X\to U$, such
that the sections $U\to U\times_C X$ induced by the $\sigma_i$ correspond
to the sections $V/\Gamma\to Y/\Gamma$.
\end{enumerate}
\end{dfn}

\begin{lem} Let $(U,Y\to V,\Gamma)$ be a chart for a family of generically
fibered surfaces $X\to C\to S$. If $S'\to S$ is an arbitrary morphism, then
$(S'\times_S U,S'\times_S Y\to S'\times_S V,\Gamma)$ (with the obvious
definitions of the various maps and of the action of $\Gamma$ on
$S'\times_S Y$ and $S'\times_S V$) is a chart for the family of generically
fibered surfaces $S'\times_S X\to S'\times_S C\to S'$. \end{lem}

\proof Conditions (1) through (4) in the definition are immediately verified
for the pullback diagram. The only point that requires a little care is to
check that $(S'\times_S 
V)/\Gamma\simeq S'\times_S U$ and $(S'\times_S Y)/\Gamma\simeq (S'\times_S
U)\times_C X$, which requires the hypothesis that $S$ be a scheme over {\bf
Q}.\endproof 

\begin{prp}
\begin{enumerate}
\item Let $(U,Y\to V,\Gamma)$ be a chart for a
family of generically fibered surfaces
$X\to C\to S$. Let $V'\subset V$ be the inverse image of $C\sm$. Then
\begin{enumerate}
\item the action of $\Gamma$ on $V'$ is free; and \item the natural
morphism $Y|_{V'}\to V'\times_S X$ is an isomorphism. \end{enumerate}
\item Furthermore, if $t_0$ is a geometric point of $S$ and $v_0$ a nodal
point of the fiber $V_{t_0}$ of $V$ over $t_0$, then \begin{enumerate}
\item the stabilizer $\Gamma'$ of $v_0$ is a cyclic group which sends each of
the branches of $V_{t_0}$ to itself.
\item If $n$ is the order of\/ $\Gamma'$, then a generator of
$\Gamma'$ acts on the tangent space of each branch by multiplication with a
primitive $n$-th root of 1. \end{enumerate}
\end{enumerate}
In particular, $V' = V\sm.$
\end{prp}

\proof The claim (1b), that the natural morphism $Y|_{V'}\to V'\times_S
X$ is an isomorphism, is a consequence of the first statement (1a).

Let $v_0$ be a geometric point of the inverse image of $U\sm$ in $V$,
$\Gamma'$ its stabilizer. By definition $\Gamma'$ acts faithfully on the
fiber $Y_0$ of $Y$ on $v_0$. The fiber $X_0$ of $X$ over the image of $v_0$
in $C$ is the quotient of $Y_0$ by $\Gamma'$; from Lemma
\ref{QuotientDropsGenus} above we get 
that $\Gamma'$ is trivial, as claimed.

For part (2) of the Proposition we observe that if the stabilizer $\Gamma'$
of $v_0$ did not preserve the branches of $V_{t_0}$ then the quotient
$V_{t_0}/\Gamma'$, which is \'etale at the point $v_0$ over the fiber
$U_{t_0}$, would be smooth over $S$ at $v_0$, so $v_0$ would be in the inverse
image of 
$U\sm$. From part (1) of the Proposition it would follow that $\Gamma'$ is
trivial, a contradiction.

So $\Gamma'$ acts on each of the two branches individually. The action on
each branch must be faithful because it is free on the complement of the
set of nodes; this means that the representation of $\Gamma'$ in each of
the tangent spaces to the branches is faithful, and this implies the final
statement.\endproof

\begin{dfn} A chart is called\/ {\em balanced} if for any nodal point of any
geometric fiber of $V$, the two roots of 1 describing the action of a generator
of the stabilizer on the tangent spaces to each branch of $V$ are inverse to
each other. 

It is easy to see that a chart is balanced if and only if it admits a
deformation to a smooth curve.\end{dfn}

\subsection{The transition scheme of two charts}
Let $X\to C\to S$ be a family of generically fibered surfaces, $(U_1,Y_1\to
V_1,\Gamma_1)$, $(U_2,Y_2\to V_2,\Gamma_2)$ two charts; call ${\rm
pr}_i\colon V_1\times_C V_2\to V_i$ the $i^{\rm th}$ projection. Consider
the scheme
$$
I = {\mathop{\isom}\limits_{V_1\times_C V_2}}({\rm pr}_1^*Y_1,{\rm pr}_2^*Y_2)
$$
over $V_1\times_C V_2$ representing the functor of isomorphisms of the two
families ${\rm pr}_1^*Y_1$ and ${\rm pr}_2^*Y_2$. There is a section of $I$
over the inverse image $\widetilde V$ of
$C\sm$ in $V_1\times_C V_2$ which corresponds to the isomorphism ${\rm
pr}_1^*Y_1\rest{\widetilde V}\simeq {\rm pr}_2^*Y_2\rest{\widetilde V}$
coming from the fact that both ${\rm pr}_1^*Y_1$ and ${\rm pr}_2^*Y_2$ are
pullbacks to $\widetilde V$ of the restriction of $X$ to $C\sm$. We will
call the scheme-theoretic closure $R$ of this section in $I$ the {\em
transition  scheme\/} from $(U_1,Y_1\to V_1,\Gamma_1)$ to $(U_2,Y_2\to
V_2,\Gamma_2)$; it comes equipped with two projections $R\to V_1$ and
$R\to V_2$. There is also an action of
$\Gamma_1\times\Gamma_2$ on $I$, defined as follows. Let
$(\gamma_1,\gamma_2)\in \Gamma_1\times\Gamma_2$, and $\phi\colon {\rm
pr}_1^*Y_1\simeq {\rm
pr}_2^*Y_2$ an isomorphism over $V_1\times_C V_2$; then define
$(\gamma_1,\gamma_2)\cdot \phi = \gamma_2\circ \phi\circ
\gamma_1^{-1}$.\footnote{Angelo - I changed the order of the gammas here, I
hope correctly}
This action of $\Gamma_1\times\Gamma_2$ on $I$ is compatible with the
action of $\Gamma_1\times\Gamma_2$ on $V_1\times_C V_2$, and leaves $R$
invariant. It follows from the definition of an essential action that the
action of $\Gamma_1 = \Gamma_1\times \{1\}$ and $\Gamma_2 =
\{1\}\times\Gamma_2$ on $I$ is free.

\begin{dfn} Two charts $(U_1,Y_1\to
V_1,\Gamma_1)$ and $(U_2,Y_2\to V_2,\Gamma_2)$ are\/ {\rm compatible} if
their transition scheme $R$ is \'etale over $V_1$ and $ V_2$.
\end{dfn}

Let us analyze this definition. Start from two charts $(U_1,Y_1\to
V_1,\Gamma_1)$ and $(U_2,Y_2\to V_2,\Gamma_2)$. Fix two geometric points
$v_1\colon \spec \Omega
\to V_1$ and $v_2\colon \spec \Omega \to V_2$ mapping to the same geometric
point $v_0\colon \spec
\Omega\to C$, and call $\Gamma'_i\subseteq \Gamma_i$ the stabilizer of
$v_i$. Also call $ V_1\sh$,
$ V_2\sh$ and
$C\sh$ the spectra of the strict henselizations of $ V_1 $, $ V_2 $ and $C$
at the points $v_1,v_2$ and $v_0$ respectively.
The action of $\Gamma_i$ on $V_i$ induces an action of $\Gamma'_i$ on $
V_i\sh$. Also call $Y_i\sh$ the pullback of $Y_i$ to $ V_i\sh$; there is an
action of $\Gamma'_i$ on $Y_i\sh$ compatible with the action of $\Gamma'_i$
on $ V_i$.

\begin{prp} The two charts are compatible if and only if for any pair of
geometric points $v_1$ and $v_2$ as above there exist an isomorphism of
groups $\eta\colon \Gamma'_1\simeq \Gamma'_2$ and two compatible
$\eta$-equivariant isomorphisms
$\phi\colon V_1\sh\simeq V_2\sh$ and
$\psi\colon Y_1\sh\to Y_2\sh$ of schemes over
$C\sh$.
\end{prp}

\proof
Consider the spectrum $(V_1\times_C V_2)\sh$ of the strict henselization of
$V_1\times_C V_2$ at the point $(v_1,v_2)\colon \spec\Omega\to V_1\times_C
V_2$, and call $R\sh$ the pullback of $R$ to $(V_1\times_C V_2)\sh$. Assume
that the two charts are compatible. The action of $\Gamma_1\times\Gamma_2$
on $I$ described above induces an action of
$\Gamma'_1\times\Gamma'_2$ on $R\sh$, compatible with the action of
$\Gamma'_1\times\Gamma'_2$ on
$(V_1\times_C V_2)\sh$. The action of $\Gamma'_1 = \Gamma'_1\times\{1\}$ on
the inverse image of $C\sm$ in $R\sh$ is free, and its quotient is the
inverse image of $C\sm$ in $V_2\sh$; but $R\sh$ is finite and \'etale over
$V_2\sh$, so the action of $\Gamma'_1$
on all of $R\sh$ is free, and $R\sh/\Gamma'_1 = V_2$. Analogously the
action of $\Gamma'_2$ on $R\sh$ is free, and $R\sh/\Gamma'_2 = V_1$.

Now, each of the
connected components of $R\sh$ maps isomorphically onto both $V_1$ and
$V_2$, because $V_i$ is
the spectrum of a strictly henselian ring and the projection $R\sh\to V_i$
is \'etale; this implies in particular that the order of $\Gamma_1$ is the
same as the number $n$ of connected components, and likewise for
$\Gamma_2$. Fix one of these components, call it
$R_0\sh$; then we get isomorphisms $R_0\sh\simeq V_i$, which yield an
isomorphism $\phi\colon V_1\simeq V_2$.

Call $\Gamma'$ the stabilizer of the
component $R_0\sh$ inside $\Gamma'_1\times\Gamma'_2$; the order of
$\Gamma'$ is at least
$|\Gamma'_1\times\Gamma'_2|/n = n^2/n = n$. But the action of $\Gamma'_2$
on $R\sh$ is free, and so $\Gamma'\cap \Gamma_2 = \{1\}$; this implies that
the order of $\Gamma'$ is $n$, and the projection $\Gamma'\to \Gamma_1$ is
an isomorphism. Likewise the projection $\Gamma'\to \Gamma_2$ is an
isomorphism, so from these we get an isomorphism $\eta\colon \Gamma_1\to
\Gamma_2$, and it is easy to check that the isomorphism of schemes
$\phi\colon V_1\simeq V_2$ is $\eta$-equivariant.

There is also an isomorphism of the
pullbacks of $Y_1\sh$ and $Y_2\sh$ to $R_0\sh$, coming from the natural
morphism $R_0\sh\to I$, which induces an isomorphism $\psi\colon Y_1\sh\to
Y_2\sh$. This
isomorphism is compatible with $\phi$, and is it also $\eta$-equivariant.

Let us prove the converse. Suppose that there exist $\eta$, $\phi$ and
$\psi$ as above. Then there is a morphism $\sigma\colon
V_1\sh\times\Gamma'_1\to I$ which sends a point $(v_1,\gamma_1)$ of
$V_1\sh\times\Gamma'_1$ into the point of $I$ lying over the point $(v_1,
\phi \gamma_1v_1) = (v_1,\eta(\gamma_1)\phi v_1)$ corresponding to the
isomorphism $\gamma_1\psi$ of the fiber of $Y_1$ over $v_1$ with the fiber
of $Y_2$ over $\phi \gamma_1v_1$. The morphism $\sigma$ is an isomorphism
of $V_1\sh\times\Gamma'_1$ with $R\sh$ in the inverse image of $C\sm$; it
also follows from the fact that the action of $\Gamma'$ on $Y_1\to V_1$ is
essential that $\sigma$ is injective. Since the inverse image of $C\sm$ is
scheme-theoretically dense in $R\sh$ and $V_1\sh\times\Gamma_1$ is
unramified over $V_1$ we see that $\sigma$ is an isomorphism; it follows
that $R\sh$ is \'etale over $V_1\sh$; analogously it is \'etale over
$V_2\sh$. So $R$ is \'etale over $V_1$ and $V_2$ at the points $v_1$ and
$v_2$; since this holds for all $v_1$ and $v_2$ mapping to the same point
of $C$ the conclusion follows.\endproof

In \ref{example1} and \ref{example2} below, we give two examples of
incompatible charts on the same coarse fibered surface.

\subsection{The product chart}
Given two compatible charts $(U_1,Y_1\to V_1, \Gamma_1)$ and $(U_2,Y_2\to
V_2, \Gamma_2)$, the graph $Y\subseteq (Y_1\times_{V_1} R)\times_R
(R\times_{V_2} Y_2)$ over $R$ of the canonical isomorphism of the two
families $Y_1\times_{V_1} R$ and $R\times_{V_2} Y_2$ is invariant under the
action of $\Gamma_1\times\Gamma_2$, and it is a family of pointed curves on
$R$. Then
$$( U_1\times_C U_2, Y\to R,
\Gamma_1\times\Gamma_2)$$
is a chart, called {\it the product chart}. It is compatible with both the
original charts.

Compatibility of charts is stable under base change:

\begin{prp} Let $(U_1,Y_1\to V_1,\Gamma_1)$ and $(U_2,Y_2\to V_2,\Gamma_2)$
be two compatible charts for a family of generically fibered surfaces $X\to
C\to S$. If $S'\to S$ is an arbitrary morphism, then
$$(S'\times_S
U_1,S'\times_S Y_1\to S'\times_S V_1,\Gamma_1)$$ and
$$( S'\times_S U_2,S'\times_S Y_2\to S'\times_S V_1,\Gamma_2)$$
are compatible charts for the the family $S'\times_S X\to S'\times_S C\to S'$.
\end{prp}

This is easy.

We now come to the definition of our basic object:

\begin{dfn} A\/ {\em family of fibered surfaces} $$
\begin{array}{c}\X\\ \down \\ \C \\ \down \\ S\end{array} $$
is a
family of generically fibered surfaces $X\to C\to S$ such that $C\to S$ is
proper, together with a collection $\{(U_\alpha, Y_\alpha\to
V_\alpha,\Gamma_\alpha)\}$ of mutually compatible charts, such that the
images of the $U_\alpha$ cover $C$.

Such a collection of charts is called an\/ {\em atlas}.

A family of fibered surfaces is called\/ {\em balanced} if each chart in its
atlas is balanced.

The family of generically fibered surfaces $X\to C\to S$ supporting the
family of fibered surfaces $\X\to \C\to S$ will be called a family of {\em
coarse fibered surfaces}.
\end{dfn}

\begin{lem} If two charts for a family of fibered surfaces are compatible
with all the charts in an atlas, they are mutually compatible.

Furthermore, if the family is balanced, then any chart which is compatible with
every chart of the atlas is balanced. \end{lem}

\smallbreak
The proof is straightforward.

\begin{rem} The lemma above allows us to define a family of fibered surfaces
using a 
maximal atlas, if we want.
\end{rem}

\begin{dfn}\label{def-morphism} A {\em morphism of fibered surfaces} $\X\to
\C\to S$ to $\X'\to \C'\to S'$   is a {\em cartesian} diagram of coarse fibered
surfaces $$ 
\begin{array}{ccc}X&\longrightarrow &X'\\ \down&&\down\\
C&\longrightarrow &C'\\
\down&&\down\\
S&\longrightarrow &S',\\
\end{array}$$
such that the pullback of the charts of an atlas of $\X'\to \C'\to S$ are all
compatible with the atlas of $\X\to \C\to S$.

The
composition of morphisms is the obvious one. \end{dfn}

We will soon reinterpret this definition of a morphism.

\subsection{\bf Fibered surfaces as stacks.} \label{fs-as-stacks}

Consider a family of fibered surfaces
$\X\to \C\to S$ with an atlas $\{(U_\alpha, Y_\alpha\to
V_\alpha,\Gamma_\alpha)\}$. For each pair of indices $(\alpha,\beta)$ let
$R_{\alpha\beta}$ be the transition scheme from $(U_\alpha, Y_\alpha\to
V_\alpha, \Gamma_\alpha)$ to $(U_\beta, Y_\beta\to V_\beta,\Gamma_\beta)$. Let
$V$ be the disjoint union of the $V_\alpha$, and let $R$ be the disjoint union
of  the $R_{\alpha,\beta}$. These have the following structure: 
\begin{itemize}
\item there are two projection $R\to V$, which are
\'etale;
\item there is a natural diagonal morphism $V\to R$ which sends each
$V_\alpha$ to $R_{\alpha\alpha}$; and
\item there is a product $R\times_V R\to R$, sending each
$R_{\alpha\beta} \times_{V_\beta} R_{\beta\gamma}$ to $R_{\alpha\gamma}$
via composition of isomorphisms
\end{itemize}
 These various maps give $R\double V$ the
structure of a groupoid. Since the diagonal map $R\to V\times_C V$
is unramified the groupoid defines a quotient Deligne-Mumford stack. In a
slight abuse of notation, we call this stack $\C$ as well.

A similar groupoid structure can be formed using $Y_\alpha$ and their pullback
to $R_{\alpha\beta}$, endowing $\X$ with the structure of a Deligne-Mumford
stack as well. Note that $\X \to \C$ is representable, and the stack $\X$ is a
family of stable $\nu$-pointed curves of genus $\gamma$ over the stack $\C$.

Let us list some properties of $\X$ and $\C$:
\begin{itemize}
\item There is a moduli morphism $\C \to \M$ associated to the family $\X \to
\C$. 

\item
The stack $\C$ is a proper nodal stack over $S$, namely, it has an \'etale
cover, given by the schemes $V$ underlying the charts, which are nodal over S.

\item
Over the inverse image of $C\sm$ the scheme $R$ is isomorphic to the fibered
product $V\times_C V$, so the quotient stack  $\C\sm$ coincides with $C\sm$. A 
similar statement holds for $\X$.

\item
Since $U_\alpha$ is the schematic quotient of $V_\alpha$ by the action of
$\Gamma_\alpha$, and $U\to C$ is \'etale, it is immediate that $C$ is the
schematic quotient of the groupoid $R\double V$, in other words, the stack
morphism  $\C \to C$ exhibits $C$ as the coarse moduli scheme of
$\C$. Similarly, $X$ is the coarse moduli scheme of $\X$.
\end{itemize}

As mentioned above,  a fibered surface being {\em balanced} is tantamount to
the existence of local smoothing of the fibered surface.  It is interesting to
interpret the requirement that the action of $\Gamma$ is {\em essential} within
the language of stacks.

We claim that  the action being essential is equivalent to the condition that
the moduli morphism $\C \to \M$ be {\em representable}. This follows from the
definition of essential action, using the following well known lemma:

\begin{lem}
Let $g:\G\to\F$ be a morphism of Deligne Mumford stacks. The following two
conditions are equivalent:
\begin{enumerate}
\item The morphism  $g:\G\to\F$ is representable.
\item For any algebraically closed field $k$ and any $\xi\in \G(k)$, the
natural group homomorphism $\aut (\xi) \to \aut(g(\xi))$ is a monomorphism.
\end{enumerate}
\end{lem}
\proof By definition, the morphism $g:\G\to\F$  is representable if and only if
the following condition holds: 
\begin{itemize}
\item For any algebraic space $X$ and any morphism $f:X \to F$, the stack
$\Y:=\G \times_\F X$ is equivalent to an algebraic space.
\end{itemize}
For fixed $f:X \to F$ latter condition is equivalent to the following:
\begin{itemize}
\item The diagonal morphism $\Delta:\Y \to \Y\times \Y$ is a monomorphism.
\end{itemize}
By definition, this means:
\begin{itemize}
\item Given an algebraically closed field $k$, and an element $q\in  \Y\times
\Y(k)$ and two elements $p_i\in \Y(k)$, with isomorphisms $\beta_i:\Delta(p_i)
\to q$, there exists a unique isomorphism $\phi:p_1 \to p_2$ such that
$\beta_1= \beta_2\circ\Delta(\phi)$. 
\end{itemize}
We can write $$q = \left((\g_1,x_1, \alpha_1),(\g_1,x_1, \alpha_1)\right)$$
where $\g_i\in \G(k), x_i\in X(k)$, and $ \alpha_i:g(\g_i) \stackrel{\sim}{\to}
f(x_i)$. Similarly $p_i=(\g'_i,x'_i, \alpha'_i)$ as above. The existence of
$\beta_i$ implies $x_1 = x_2 = x'_1 = x'_2$. Also, composing with  the given
isomorphisms 
we may reduce to the case that, in fact,  the same is true for the $\g_i,
\g_i', \alpha_i$ and $\alpha_i'$. Thus the condition above is equivalent to
saying that for any 
$\g\in G(k)$ there is a unique $\beta\in \aut(\g)$ such that $g(\beta) =
id_{g(\g)}$, which is what we wanted. \qed

\subsection{A stack-theoretic formulation of the category of fibered surfaces}

We will now formulate  stack-theoretic data similar to that obtained above, and
 then compare them to  the data of a fibered surface.

Let $S$ be a scheme over {\bf Q}. Consider a proper
Deligne--Mumford stack ${\cal C}\to S$, such that its fibers
are purely one-dimensional and geometrically connected, with
nodal singularities. Call $C$ the moduli space of ${\cal C}$;
this automatically exists as an algebraic space.

\begin{prp} The morphism $C \to S$ is a flat
family of nodal curves.
\end{prp}

\proof First of all let us show that $C$ is flat over $S$. We
may assume that $S$ is affine; call $R$ its ring. Fix a
geometric point $c_0 \to C$, and call
$C\sh$ the strict henselization of $C$ at $c_0$. Let $U$ be an
\'etale cover of ${\cal C}$, $u_0$ a geometric point of $U$
lying over $c_0$, and call $U\sh$ the strict henselization of
$U$ at $u_0$. If $\Gamma$ is the automorphism group of the
object of ${\cal C}$ corresponding to $u_0$, then $\Gamma$
acts on $U\sh$, and $C\sh$ is the quotient $U\sh/ \Gamma$.
Because the schemes are defined over {\bf Q}, the ring of
$C\sh$ is a direct summand of the ring of $Y\sh$, as an
$R$-module, so it is flat over $R$.

The fact that the fibers are nodal follows from the fact
that, over an algebraically closed field, the quotient of a
nodal curve by a group action is again a nodal curve.\endproof

\begin{dfn} A stack-like family of
curves of genus $g$ over a scheme $S$ consists of a
proper stack ${\cal C}\to S$, whose fibers are purely
one-dimensional and geometrically connected, with nodal
singularities, such that:

\begin{enumerate}
\item the fibers of the morphism $C \to S$, where $C$ is
the moduli space of ${\cal C}$, have genus
$g$;

\item over the inverse image of the smooth locus $C\sm$
of the map $C \to X$, the projection ${\cal C} \to C$ is an
isomorphism.

\end{enumerate}

If ${\cal C} \to S$ and ${\cal C}' \to S'$ are stack-like
families of curves of genus $g$, a
morphism $F \colon {\cal C}\to{\cal C}'$ consists of a
cartesian diagram
$$
\begin{array}{ccc}
{\cal C}&\mapright F &{\cal C}'\\
\down&& \down\\
S& \mapright f &S'.
\end{array}
$$

\end{dfn}

The composition of morphisms of stack-like families of 
curves is defined in the obvious way. In this way stack-like
families of  curves of genus $g$ form a
2-category. The 2-arrows being defined in the usual way.

\begin{prp} The 2-category of stack-like families of
 curves of genus $g$ is equivalent to a
1-category.
\end{prp}

\proof This is the
same as saying that a 1-arrow in the category can not have
nontrivial automorphisms. The point here is that the stack
${\cal C}$ has an open subscheme which is
stack-theoritically dense in ${\cal C}$, which is sufficient by the following
lemma.

\begin{lem} Let $F \colon  {\cal X} \to {\cal Y}$ be a representable
morphism of Deligne--Mumford stacks over a scheme $S$. Assume
that there exists an open representable substack (i.e. an algebraic space) $U
\subseteq {\cal X}$ and a dense open representable substack $V
\subseteq {\cal Y}$  such that $F$ maps $U$ into $V$. Further assume that the
diagonal 
${\cal Y} \to {\cal Y}\times_S {\cal Y}$ is separated. Then
any automorphism of $F$ is trivial.
\end{lem}

\proof 
First note that the lemma holds if $ {\cal X} = X$ is an algebraic space: if we
denote by $\eta$ the 
object of ${\cal Y}$ over $X$ corresponding to $F$, then 
the representability of the diagonal ${\cal Y} \to {\cal Y}\times_S {\cal Y}$
implies that the isomorphism scheme $\isom_X(\eta,\eta)\to X$ is
separated. Since $V$ is repesentable, we have that $\isom_U(\eta|_U,\eta|_U)\to
U$ is an isomorphism. Thus the unique section over the given open set $U$ has
at most one extension to $X$, which gives the assertion in this case. We
will now  reduce the general case to this one by descent. We start with some
observations.

Let $\alpha$ be an
automorphism of $F$; for each object $\xi$ of ${\cal X}$ over
a scheme $X$ we are given an automorphism $\alpha_\xi$ of
$F(\xi)$ over $X$, satisfying the usual condition for being a
natural transformation. We are going to need the following
two facts. First of all, if $\xi \to \zeta$ is a morphism in
${\cal X}$, then it follows from the fact that ${\cal X}$ is
a category fibered in groupoids that if $\alpha_\zeta$ is
trivial then also $\alpha_\xi$ is trivial.
Also, if $X' \to X$ is an
\'etale surjective map of schemes and $\xi'$ is the pullback
of $\xi$ to $X'$, then it follows from the definition of a
stack that if $\alpha_{\xi'}$ is trivial then also
$\alpha_\xi$ is trivial.

 Let $X \to {\cal X}$ be an \'etale
cover, and call $\xi$ the corresponding object of ${\cal X}$
over $X$. The restriction of $\alpha_\xi$ to the open
subscheme $\tilde{U} = U \times_{\cal X} X \subseteq X$ is trivial,
and $\tilde{U}$ is scheme-theoretically dense in $X$; applying this lemma in
the case of algebraic spaces it follows that
$\alpha_\xi$ is trivial.

Take an arbitrary object $\tau$ of ${\cal X}$ over a scheme
$T$; then there is an \'etale cover $T' \to T$ such that the
pullback $\tau'$ of $\tau$ to $T'$ admits a morphism $\tau
\to \xi$. By applying the two facts mentioned above we get
the result.\endproof

\begin{dfn} A stack-like fibered surface $f \colon
{\cal C}\to\M$ over a scheme $S$ is a stack-like family ${\cal
C}\to S$ of
 curves of genus $g$ over $S$, with projective coarse moduli space $C\to S$,
together 
with a representable morphism ${\cal C} \to \M$.
\end{dfn}

As pointed out above, stack-like fibered surfaces form a 2-category.  A 1-arrow
$\Phi= (\phi, \alpha_\Phi)$ from 
$f \colon {\cal C}\to\M$ to $f' \colon {\cal C}'\to\M$
consists of a morphism of stack-like families of curves $\phi
\colon {\cal C}\to {\cal C}'$, together with an isomorphism $\alpha_\Phi$
of the composition of $f' \circ \phi$ with
$f$; the compostion of 1-arrows is defined in the obvious
way. 
A 2-arrow from a 1-arrow $\Phi$ from $f
\colon {\cal C}\to\M$ to $f' \colon {\cal C}'\to\M$ to
another arrow
$\Psi$ with the same domain and codomain consists of an
isomorphism $\beta$ of the two functors $\phi$ and $\psi$
, with the usual compatibility condition of the isomorphism
$\beta$ with $\alpha_\Phi$ and $\alpha_\Psi$.

Because of the Proposition, this category is equivalent to a
nice unproblematic 1-category, which we can think of as the
category of stack-like fibered surfaces, the arrows being
given by isomorphism classes of morphisms.

We can now compare the two categories we have defined.

\begin{th} The category of stack-like fibered
surfaces over schemes over {\bf Q} is equivalent to the
category of fibered surfaces.
\end{th}

\proof 
The construction in \ref{fs-as-stacks} shows that a fibered surface gives rise
to a stack-like fibered surface in a functorial way.






To go from stack-like fibered surfaces to honest fibered
surfaces, let us take a stack-like fibered surface ${\cal C}
\to \M$, and then take as a family of pointed nodal curves
its coarse moduli space $C \to S$. The charts are given by
\'etale morphism from schemes to ${\cal C}$. We leave the treatment of
morphisms to the reader.\qed

%
%

From now on we will use fibered surfaces and stack-like fibered surfaces
interchangeably. 

\subsection{Stable fibered surfaces}

We first define a moduli morphism on a coarse fibered surface:

\begin{lem} Let $\X\to \C\to S$ be a family of fibered surfaces. Then the
morphism $C\sm\to \m$ induced by the restriction of $X\to C$ to $C\sm$
extends uniquely to a morphism $C\to \m$.
\end{lem}

\proof The unicity is clear from the fact that $\m$ is separated and $C\sm$
is scheme-theoretically dense in $C$. To prove the existence of an
extension is a local question in the \'etale topology; but if $\{(
U_\alpha,Y_\alpha\to V_\alpha,\Gamma_\alpha)\}$ is an atlas then  the families
$Y_\alpha\to V_\alpha$ 
induce morphisms $V_\alpha\to \m$, which are $\Gamma_\alpha$-equivariant,
yielding 
morphisms $U_\alpha\to \m$. These morphisms are extensions of the restriction
to the 
$U_\alpha$ of the morphism $C\sm\to \m$. Therefore they descend to
$C$.\endproof 

We can use this lemma to finally define stable fibered surfaces:

\begin{dfn} A family of fibered surfaces $\X\to \C\to S$ is {\em stable} if
the associated morphism $C\to \m$ is Kontsevich stable.
\end{dfn}

\subsection{Natural line bundles on a coarse fibered surface}

\begin{prp}\label{rel-ample-sheaf} To each family of fibered surfaces $\X\to
\C\to 
S$ it is possible 
to assign a canonically defined line bundle ${\cal L}_X$
on the coarse fibered surface $X$, which is relatively ample along the map $X
\to 
C$. This line bundle satisfies:

\begin{enumerate} 
\item For any morphism of fibered surfaces
$$
\begin{array}{ccc}
\X'&\mapright\phi &\X\\
\down&&\down\\
\C'&\longrightarrow &\C\\
\down&&\down\\
S'&\longrightarrow &S\\
\end{array}
$$
there is an isomorphism
of line bundles $\sigma_\phi\colon {\cal L}_{X'}\simeq \phi^*{\cal L}_X$.

\item These isomorphisms satisfy the cocycle condition, in the sense that
$\sigma_{{\rm id}_X} = {\rm id}_{{\cal L}_X}$ for all families of fibered
surfaces $X\to S$, and if we have two morphisms $$
\begin{array}{ccccc}
\X''&\mapright\psi&\X'&\mapright\phi &\X\\
\down&&\down&&\down\\
\C''&\longrightarrow& \C'&\longrightarrow& \C\\ \down&&\down&&\down\\
S''&\longrightarrow & S'&\longrightarrow &S\\ \end{array}
$$
then
$$
\sigma_{\phi\psi} =
\psi^*\sigma_\phi\circ\sigma_\psi \colon {\cal L}_{X''}\longrightarrow
(\phi\psi)^*{\cal L}_X = \psi^*\phi^*{\cal L}_X
$$
\end{enumerate}
\end{prp}

\proof Fix an integer $N$ which is divisible by the order of the
automorphism group of any $\nu$-pointed stable curve of genus $\gamma$. The
construction below will depend on the choice of the integer $N$. For
each index $\alpha$ consider the line bundle $\omega_\alpha =
\omega_{Y_\alpha/V_\alpha}(\Delta_\alpha)$, the relative dualizing sheaf of
$Y_\alpha$ over $V_\alpha$ twisted by the divisor $\Delta_\alpha$ of marked
points. This line bundle has a natural action of $\Gamma_\alpha$, and the
stabilizers of all the geometric points act trivially on the fiber of
$\omega_\alpha^{\otimes N}$; therefore $\omega_\alpha^{\otimes N}$ descends
to a line bundle ${\cal L}_\alpha$ on $X\times_C U_\alpha$.

Given two indices $\alpha$ and $\beta$, let $R_{\alpha\beta}$ be the
transition scheme of two corresponding charts and $Y_{\alpha\beta}\to
R_{\alpha\beta}$ the family of pointed stable curves constructed above; there
is a canonical 
isomorphism between the pullback of $\omega_\alpha$ to $Y_{\alpha\beta}$
and the analogous bundle $\omega_{\alpha\beta}$ on the family
$Y_{\alpha\beta}$, so we can compose these to get isomorphisms between the
pullbacks of $\omega_\alpha$ and $\omega_\beta$. These isomorphisms satisfy
the cocycle condition and give the collection of the $\omega_\alpha$ the
structure of a line bundle on the groupoid $R\double V$.

Passing to the $N$-th power, the isomorphisms above also descend to
isomorphisms 
between the pullbacks of ${\cal L}_\alpha$ and ${\cal L}_\beta$ on the
fiber product $X\times_C (U_\alpha\times_C U_\beta)$ satisfying the cocycle
condition, so the ${\cal L}_\alpha$ are the pullback of a well defined line
bundle ${\cal L}_X$, which is the one we want.\endproof

\begin{prp}\label{ample-sheaf} To each family of stable fibered surfaces $\X\to
\C\to S$ it is possible  
to assign a canonically defined line bundle ${\cal A}_X$
on the coarse fibered surface $X$, which is relatively ample along the map $X
\to 
S$. This line bundle satisfies the properties listed in Proposition
\ref{rel-ample-sheaf}.
\end{prp}

\proof Fix $N$ as in the proof of Proposition 
\ref{rel-ample-sheaf}. Fix an ample line bundle $H_0$ on $\m$, and let $H =
H_0^{\otimes 3}$. Let $\X\to
\C\to S$ be a stable fibered surface, and let $m_C:C\to \m$ be the associated
moduli morphism. Consider the line bundle $A_C:=(\omega_{C/S}\otimes
m_C^*H)^{\otimes N}$ on $C$. By the theory of stable maps, the line bundle
$A_C$ is ample relatively to the morphism $C \to S$. It clearly satisfies the
invariance condition in the proposition. Also, by Koll\'ar's semipositivity
lemma (see \cite{kollar-projectivity}) the line bundle ${\cal L}_X$ constructed
above 
is nef. Set $ {\cal A}_X = {\cal L}_X \otimes
f^*A_C$. By Kleiman's criterion for ampleness we clearly have that  $ {\cal
A}_X$ is ample. The invariance
conditions follow by construction.

\begin{rem} It is important to note that both $ {\cal L}_X $ and ${\cal A}_X$
are invariants of the {\em coarse} fibered surfaces, independent of the given
atlas. This is because, on the open subscheme $U\subset X$ where $X \to S$ is
Gorenstein, these coincide with $\omega_{U/C}^N$ and $\omega_{U/S}^N$. Since
$U$ has complement of codimension $\geq 2$, the extensions to $X$ are unique.
\end{rem}

\section{The stack of fibered surfaces}

Consider the stack $\M$ of stable $\nu$-pointed curves of genus $\gamma$ over
{\bf Q}, and 
the associated moduli space $\m$; choose an ample line bundle ${\cal H}$ on
$\m$. Fix two nonnegative integers $g$ and
$d$. 

We define a category $\PP$, fibered over the category
$\operatorname{Sch}/\bfq$ of schemes  over {\bf Q},
as follows. The objects are 
stable families $\X\to \C\to S$ such that for the associated morphism
$f\colon C\to \m$, the
degree of the line bundle $f^*{\cal H}$ on each fiber of $C$ over $S$ is
$d$. The arrows are morphism of fibered surfaces.
We also have a subcategory  $\PP^{balanced}$ of stable {\em balanced}
families. 

There is an obvious morphism from $\PP$ to the stack $\K$ of Kon\-tse\-vich
- stable maps of genus $g$ and degree $d$ into $\m$ which sends
each stable family of fibered surfaces $\X\to \C\to S$ to the associated
morphism $C\to \m$. 

Our main result is:

\begin{th}\label{main-theorem}
The category $\PP$ is a complete Deligne-Mumford stack, admitting a projective 
coarse moduli space $\pp$. The subcategory  $\PP$ forms an open and closed
substack. 
\end{th}

\proof

By definition, $\PP$ is fibered by groupoids over $\operatorname{Sch}/\bfq$,
since all the morphisms of objects over a scheme $S$ are isomorphisms. 

The next step is to check the sheaf axioms, to prove that $\PP$ is indeed a
stack. First, we need to show that we have effective descent for
objects (see \cite{artin}, 1.1 (ii)). Given a scheme $S$, an
\'etale cover $\{S_\alpha\to S\}$, and fibered surfaces
${\cal X}_\alpha$ over $S_\alpha$ together with isomorphisms
between the pullbacks of ${\cal X}_\alpha$ and ${\cal
X}_\beta$ to $S_\alpha\times_S S_\beta$ satisfying the
cocycle condition, we claim that these descend to a fibered surface $\cal
X$ on $S$. The existence of the coarse fibered surface $X$ as an
algebraic space is immediate; its 
projectivity follows from the existence of the ample line 
bundles ${\cal A}_{{\cal X}_\alpha}$ (Proposition
\ref{ample-sheaf}). Now, any chart
$$\begin{diagram}
\node{Y_\alpha}\arrow{s}\arrow{e}
	\node{X_\alpha\times_{C_\alpha}U_\alpha}\arrow{s}\arrow{e} 
		\node{X_\alpha}\arrow{s} \\
\node{V_\alpha}\arrow{e}\arrow{ese}
	\node{U_\alpha}\arrow{e}\arrow{se} 
		\node{C_\alpha}\arrow{s} \\
		\node[3]{S_\alpha}
\end{diagram}
$$

immediately gives rise to a chart 
$$\begin{diagram}
\node{Y_\alpha}\arrow{s}\arrow{e}
	\node{X\times_{C}U_\alpha}\arrow{s}\arrow{e} 
		\node{X}\arrow{s} \\
\node{V_\alpha}\arrow{e}\arrow{ese}
	\node{U_\alpha}\arrow{e}\arrow{se} 
		\node{C}\arrow{s} \\
		\node[3]{S}
\end{diagram}
$$
The charts coming from one $X_\alpha$ are obviousely compatible. Compatibility
of the charts 
coming from $X_\alpha$ and $X_\beta$ follows since they are compatible
over $X_{\alpha\beta}$.

The other  sheaf axiom (see \cite{artin}, 1.1 (i)) requires the $\isom$
functor to be a sheaf in the 
\'etale topology. Since we would like to show that $\PP$ is an {\em
algebraic} stack, 
we might as well address the representability of the
$\isom$ functor at the same time. 

Take two stable families of fibered surfaces ${\cal X}_1$ and
${\cal X}_2$ over the same scheme $S$; we want to show
that the functor $\isom_S({\cal X}_1,{\cal X}_2)$ of
isomorphisms of fibered surfaces is representable by a scheme
of finite type over $S$. By Grothendieck's theorem,
$\isom_S(X_1,X_2)$ and
$\isom_S(C_1,C_2)$ are represented by schemes locally of
finite type over $S$. There is an obvious functorial embedding
$$
{\mathop{\isom}\limits_S}({\cal X}_1,{\cal X}_2)\into
{\mathop{\isom}\limits_S}(X_1,X_2) { \mathbin{\mathop{\times}\limits_S} }
{\mathop{\isom}\limits_S}(C_1,C_2); 
$$
we will prove that $\isom_S({\cal X}_1,{\cal X}_2)$ is a
closed subscheme of $\isom_S(X_1,X_2)\times_S
\isom_S(C_1,C_2)$. This shows that $\isom_S({\cal X}_1,{\cal
X}_2)$ is a scheme locally of finite type over $S$; the fact
that a stable fibered surface has only finitely many
automorphisms means that $\isom_S({\cal X}_1,{\cal X}_2)$ is
quasifinite over $S$, so it is of finite type.

First of all, there is a closed subscheme
$I\subseteq \isom_S(X_1,X_2)\times_S \isom_S(C_1,C_2)$
representing the functor of compatible pairs of isomorphisms
$X_1\simeq X_2$ and $C_1\simeq C_2$; and  $\isom_S({\cal X}_1,{\cal X}_2)$ is 
obviously 
contained in $I$. So we only need to prove that given a
scheme $S$, two fibered surfaces ${\cal X}_1$ and ${\cal
X}_2$ on $S$, and compatible isomorphisms
$f\colon X_1\simeq X_2$ and
$g\colon C_1\simeq C_2$, there exists a closed subscheme
$S'\subseteq S$ such that given a morphism $T\to S$, the
pullbacks of $f$ and $g$ to $T$ give an isomorphism of
fibered surfaces if and only if $T\to S$ factors through
$S'$. By considering the given atlas on $X_1$ and the
pullback of the atlas of $X_2$ we see that we need the
following.

\begin{lem} Let $X\to C\to S$ be a coarse fibered
surface with two atlases ${\cal A}_1$ and ${\cal A}_2$. Then
there exists a closed subscheme $S'\subseteq S$ such that
given a morphism $T\to S$, the pullbacks of ${\cal A}_1$ and
${\cal A}_2$ to $T\times_S X$ are compatible if and only if
$T\to S$ factors through $S'$.
\end{lem}

\proof We will prove a local version of the fact above.

\begin{lem} Let there be given a coarse fibered surface
$X\to C
\to S$ and two charts $(U,Y\to V, \Gamma)$ and
$(U',Y'\to V', \Gamma')$. Let there be given a geometric point
$s_0$ of
$S$; assume that the fiber of $C$ over $s_0$ has a
unique nodal point
$c_0$, and that there are two unique geometric
points $v_0\in V$ and $v'_0 \in V'$ over $c_0$. Then, after
passing to an
\'etale neighborhood of $s_0$, there exists a closed subscheme
$S'\subseteq S$ such that given a morphism $T\to S$, the
pullbacks of the two charts to $T\times_S X$ are compatible
if and only if $T\to S$ factors through $S'$.
\end{lem}

The local lemma clearly implies the global lemma. 


{\bf Proof of the local lemma.}\enspace By passing
to the fiber product $U
\times_C U'$ we may assume that $C = U = U'$. By
refining $V$ and $U$ we may assume that they are affine,
and that there exists an invariant effective Cartier divisor
$D$ on
$V$ containing the locus where the projection
$V\to S $ is not smooth, but none of the fibers of $V\to S$.

We claim that such a divisor is flat over $S$. Indeed, we may assume
that $D$ is defined by an equation $f$:
$$ 0 \to \co_V \stackrel{\cdot f}{\to}\co_V \to \co_D \to 0.$$
Since $f$ does not vanish along any component of any fiber we have that the
natural 
map $$ \co_V{\mathop{\otimes}\limits_{\co_S}}\co_s
\stackrel{\cdot f}{\to}\co_V{\mathop{\otimes}\limits_{\co_S}}\co_s $$
is injective for any point $s\in S$.

 By passing to an
\'etale neighborhood of $s_0$ we can split $D$ into a number
of connected components, so that the components containing
$v_0$ is finite over $S$; then by deleting the other
connected components we see that we may assume that $D$ is
finite over $S$.

Consider the transition scheme $R$ of the two charts. There
is a free action of the group $\Gamma'$ over $R$; set $E =
R/ \Gamma'$. The projection $E\to V$ is an isomorphism over
the smooth locus of $C$, and is an isomorphism if and only if
the projection $R\to V$ is \'etale. Take a morphism $T\to S$;
then the transition schemes $R_T$ of the pullbacks of the two
charts to $T$ is the scheme-theoretic closure of the
inverse image of the smooth locus of $C$ in $R \times_S
T$, so $R_T$ is \'etale over $V\times_S T$ if and only if the
projection $(R \times_S T)/ \Gamma'  = E\times_S T\to V
\times_ST$ has a section.

Set $S = \spec \Lambda$, $V = \spec A$, $E = \spec B$,
and call $I$ the ideal of $D$ in $A$. The coordinate ring of the
complement of $D$ inside $V$ is $A' = \cup_{i=0}^\infty
I^{-i}$, and there is a natural homomorphism $B\to A'$;
given a $\Lambda$-algebra $L$ the coordinate ring of the transition
scheme of the pullbacks of the two charts to $\spec L$ is the
image of $B \otimes_\Lambda L$ in $A' \otimes_A L$, so $(R
\times_S T)/ \Gamma'  = E\times_S T\to V
\times_ST$ has a section if and only if the image of $B
\otimes_\Lambda L$ in $A' \otimes_\Lambda L$ is equal to $A
\otimes _\Lambda L\subseteq A' \otimes_\Lambda L$. Take
a set of generators of $B$ as an $A$-algebra, and call
$e_1,\ldots,e_n$ their images inside $A'/A$; the condition
that the image of $B \otimes_\Lambda L$ be equal to $A
\otimes _\Lambda L$ is equivalent to the condition that the
image of the $e_i$ in $(A'/A)\otimes_\Lambda L$ be zero. Fix
an integer $n$ such that the $e_i$ are all contained in
$I^{-n}/A$; then the sequence
$$
0\longrightarrow I^{-n}/A\longrightarrow A'/A\longrightarrow
\cup_{i=n}^\infty I^{-i}/I^{-n}\longrightarrow 0
$$
is exact and $\cup_{i=n}^\infty I^{-i}/I^{-n} =
I^{-n}\otimes_A(A'/A)$ is flat over $\Lambda$, so the
sequence stays exact after tensoring with $L$. The conclusion
is that
$R_T$ is \'etale over $V\times_S T$ if and only if the images
of the $e_i$ in $(I^{-n}/A)\otimes_\Lambda L$ is zero. But
$I^{-n}/A$ is finite and flat as a $\Lambda$-module, so it is
projective, and can be embedded as a direct summand of a free
$\Lambda$-modules $F$. If $J \subseteq \Lambda$ is the ideal
generated by the coefficients of the $e_i$ with respect
to  a basis of $F$, it is clear that the closed subscheme $S'=
\spec (\Lambda/J)$ has the desired property. \qed

Now let us conclude the proof that
$\PP$ is an algebraic stack, by verifying Artin's axioms
(\cite{artin}, Corollary~5.2). First note that Schlessinger's
first condition, in the 
strong form (\cite{artin} (2.3) condition (S 1')), follows from standard
principles. The finiteness condition (\cite{artin} (2.5) condition (S
2)) will follow from the properness of $\C$ once we construct the
deformation space below. The  
representability of the diagonal (\cite{artin}, Corollary~5.2,
condition 1) has just been checked.

Let us verify condition~3 of \cite{artin}, Corollary~5.2 about formal
completions. 

Let $A$ be a complete local ring with residue field $k$
and maximal ideal $\mm$; set
$A_n = A/\mm^{n+1}$, $S_n = \spec A_n$. Let there be given a
fibered surface ${\cal X}_n$ on $S_n$ for each $n$, together
with isomorphisms of the pullback of ${\cal X}_n$ to $S_{n-1}$
with ${\cal X}_{n-1}$. From Grothendieck's existence
theorem and the existence of a compatible system of very
ample line bundles on the
$C_n$ we get that the formal scheme $\{C_n\}$ is
algebraizable, that is, there exists a flat family of nodal
curves $C$ over $\spec A$ together with isomorphisms $C_n
\simeq C \times_S S_n$ compatible with the embeddings
$C_{n-1} \into C_n$; this
$C$ is unique up a unique isomorphism. Likewise, we get a
unique scheme $X$ on $C$, inducing the formal scheme $\{X_n\}$
over $C$.

For each $n$ consider the associated stack ${\cal C}_n$.

\begin{lem} There is a nodal curve $C'_0$, not
necessarily geometrically connected, and a flat morphism of
$k$-schemes $C'_0 \to {\cal C}_0$, which is \'etale over an
open subset of $C_0$ containing of the nodes.
\end{lem}

\proof Assume that $k$ is algebraically closed; assume also
that we have chosen an isomorphism of the group of roots of 1
in $k$ with ${\bf Q}/ {\bf Z}$. For each nodal point $p$ of
$C$ choose a chart $(U_p,Y_p \to C_p,
\Gamma_p)$ such that $V_p$ has a single point over $p$; then
$\Gamma_p$ is a finite cyclic group. Choose a cyclic group
$\Gamma$ whose order is divisible by the order of all the
$\Gamma_p$, and an embedding $\Gamma_p\subseteq \Gamma$ for
each $p$. Let $\Delta$ be a finite set of smooth points of
$C_0$ such that the pointed curve $(C_0, \Delta)$ is stable,
and call $N$ the set of nodal points of $C_0$. Call $\Pi$
the abelianization of the algebraic fundamental group of $C_0
\setminus(N \cup \Delta)$, and for each nodal point $p$ call
$A_p$ and $B_p$ the spectra of the completions of the local
rings of the normalization of $C$ at the two branches of $C$
at $p$. The fundamental groups of $A_p$ and $B_p$ are
canonically isomorphic to the free rank one
profinite group $\widehat{\bf Z}$; call $\alpha_p$ and
$\beta_p$ the images in $\Pi$ of the canonical generators
under the canonical homomorphisms $\pi_1(A_p)\to \Pi$ and
$\pi_1(B_p)\to \Pi$. The fact that $C_0
\setminus(N \cup \Delta)$ is stable insures that the
$\alpha_p$ and $\beta_p$ are part of a basis of the free
abelian profinite group $\Pi$, so there exists a continuous
homomorphism $\Pi \to \Gamma$ sending $\alpha_p$ to a
generator of $\Gamma_p \subseteq \Gamma$ and $\beta_p$ to the
opposite generator; call $C'_0$ the associated ramified
covering of $C_0$.

Over the inverse image of $N$ the
morphism ${\cal C}_0\to C_0$ is an isomorphism, so the
morphism $\pi\colon C'_0\to C_0$ yields a morphism
$C'_0\setminus \pi^{-1}(N)
\to {\cal C}_0$. There is also an isomorphism of some
\'etale neighborhood of any point of $C'_0$ mapping to some
$p\in N$ with an \'etale neighborhood of the inverse image
of $p$ in $V_p$, and each $V_p$ comes with an \'etale
morphism $V_p \to {\cal C}_0$; this gives an extension of the
morphism $C'_0\setminus \pi^{-1}(N) \to {\cal C}_0$ over the
inverse image of $N$ to all of $C'_0$, which is \'etale
everywhere except possibly at the points of $\Delta$. This
morphism is automatically flat.

If $k$ is not algebraically closed, then consider the
pullback of the fibered surface ${\cal X}_0$ to the
algebraic closure $\overline{k}$. Given a morphism
$\overline{C}'_0
\to {\cal C}_0\times_{ \spec k} \spec \overline{k}$
satisfying the desired conditions, this will be obtained by
base change from some morphism $C'_0 \to {\cal C}_0\times_{
\spec k} \spec k'$, where
$k'$ is a finite extension of $k$. The desired morphism is
the composition of this  $C'_0 \to {\cal C}_0\times_{ \spec
k} \spec k'$ with the projection ${\cal C}_0\times_{
\spec k} \spec k' \to{\cal C}_0$.
\endproof

Now starting from $C'_0$ we can construct, by induction on
$n$, a sequence of flat nodal curves $C'_n$ on $S_n = \spec
A_n$, together with embeddings $C'_{n-1} \into C'_n$,
inducing an isomorphism of
$C'_{n-1}$ with the restriction of $C'_n$ to $S_{n-1}$, and
proper flat morphisms $\phi_n\colon C'_n\to {\cal C}_n$ such that
the restriction of $\phi_n$ is $\phi_{n-1}$; for each $n$ call
$C''_n$ the fiber product $C'_n \times_{ {\cal C}_n} C'_n$;
then the groupoid $C''_n \double C'_n$ is a flat presentation
of the stack ${\cal C}_n$. There are flat schemes
$C'$ and $C''$ over $A$, whose associated formal schemes are
precisely $\{C'_n\}$ and $\{C''_n\}$. For each $n$ there is
a family ${\cal Y}_n \to {\cal C}_n$ of $\nu$-pointed stable
curves of genus
$\gamma$ on the stack ${\cal C}_n$; by pullback we get
families of stable curves $Y'_n \to C'_n$ and $Y''_n \to
C''_n$, which yield families of pointed stable curves $Y' \to
C'$ and
$Y'' \to C''$.

The induced flat groupoid
$C''\double C'$ has a quotient Deligne--Mumford stack ${\cal
C}$ over $\spec A$; the groupoid $Y''\double Y'$ gives a flat
family of pointed stable curves ${\cal X} \to {\cal C}$, or,
in other words, a morphism ${\cal C}\to \M$. Let us check
that the moduli space of ${\cal X}$ is $X$. In fact, the
moduli space of ${\cal X}_n$ is $X_n$, and if we are given a
morphism ${\cal X} \to Z$, where $Z$ is an algebraic space
over $S$, we get a compatible sequence of morphisms ${\cal
X}_n \to Z$, and from this a compatible sequence of morphisms
$X_n \to Z$. The graphs of these morphisms give a formal
subscheme of $X \times_S Z$, which corresponds to a subscheme
of $X \times_S Z$, the graph of a morphism $X \to Z$. The
composition of this with the morphism ${\cal X} \to X$ is the
given morphism ${\cal X} \to Z$.

Call $V$ the locus of $C'$ where the morphism $C' \to {\cal
C}$ is \'etale, $U = V/
\Gamma$, and call $Y$ the pullback of the family
${\cal X} \to {\cal C}$ to $V$. Then $(U, Y \to V, \Gamma)$ is
a chart for the coarse fibered surface $X \to C$; it is easy
to check that it is balanced, using the fact that the central
fiber is balanced.

Now we have to define an obstruction theory ${\cal O}$ for
the stack $\PP$. This can be done, using the results of
\cite{illusie}. A direct more elementary approach is also
possible, along the lines of \cite{vistolidef}.

Consider a reduced ring $A_0$, and a surjective homomorphism
$A' \to A$ of infinitesimal extensions of $A_0$, whose kernel
$M$ is an $A_0$-module. Consider a fibered surface ${\cal X}
\to {\cal C}$ over $\spec A$; to give a lifting of ${\cal X}$
to $\spec A'$ is equivalent to giving a flat lifting of the
morphism ${\cal C}\to \M$ to $\spec A'$. 

This induces a representable morphism of stacks $f \colon
{\cal C} \to
\M$; we consider the cotangent complex ${\cal L} = {\rm L}_{
{\cal C}/\M\times_ {\spec{\bf Q}}\spec A}$ over $A$. Because
$\M\times_ {\spec{\bf Q}}\spec A$ is smooth over $A$ and
${\cal C}$ is a local complete intersection, it follows
that the complex ${\cal L}$ is quasi-isomorphic to the
two-step complex having $f^* \Omega_{\M/ {\bf Q}}$ in degree
$-1$ and $\Omega_{{\cal C}/A_0}$ in degree 0, the differential
$f^* \Omega_{\M/ {\bf Q}}\to \Omega_{{\cal C}/A_0}$ being the
obvious homomorphism. In what follows we will
consider all sheaves over ${\cal C}$ as sheaves over the
\'etale site of ${\cal C}$, and the various cohomology and
$\ext$ groups will be considered with respect to this site.

It follows from
\cite{illusie}, Proposition~2.1.2.3., that, in the notation of
\cite{artin}, the
$A_0$-module ${\rm D}_{\cal X}(M)$ is $\ext_{{\cal
O}_{\cal C}}^1({\cal L},{\cal
O}_{\cal C} \otimes_A M)$, and that we get an
obstruction theory by setting

$$
{\cal O}_ {\cal X}(M) = \ext_{{\cal
O}_{\cal C}}^2({\cal L},{\cal
O}_{\cal C} \otimes_A M).
$$
We have to check that these $A_0$-modules are finite, and
that they satisfy the conditions of \cite{artin} (4.1).

If ${\cal K}$ is the complex
$\cursrhom_{{\cal O}_{{\cal C}}}({\cal L},{\cal
O}_{{\cal C}_0})$ of sheaves over ${\cal C}_0$, which is a
bounded complex with coherent cohomology on ${\cal C}_0$, then
we have
$$
\cursrhom_{{\cal
O}_{{\cal C}}}({\cal L},{\cal
O}_{{\cal C}_0} \otimes_{A_0}M) = {\cal K} \tototi_{A_0} M
$$
and
${\rm D}_{\cal X}(M)$ and ${\cal O}_ {\cal X}(M)$ are
respectively the first and second hypercohomology group of the
total tensor product complex ${\cal K}\tototi_A M$.

Consider
the projection $\pi\colon {\cal C} \to C$. The functor
$\pi_*$ carries sheaves of ${\cal O}_{\cal C}$-module to
sheaves of ${\cal O}_C$-modules.

\begin{lem} The functor $\pi_*$ carries
quasicoherent sheaves to quasicoherent sheaves, coherent
sheaves to coherent sheaves, and is exact.
\end{lem}

\proof The question is local in the \'etale topology on $C_0$,
so we may pass to a chart, and assume that ${\cal C}$ is of
the form $[V/\Gamma]$, where $V$ is a scheme and $\Gamma$ a
finite group. Then sheaves on ${\cal C}$ correspond to
equivariant sheaves on $V$, and the statement is
standard.\endproof

The lemma implies that the hypercohomology of the complex
${\cal K}\tototi_A M$ is isomorphic to the hypercohomology of
the comples of sheaves $\pi_*({\cal K}\tototi_A M) =
(\pi_*{\cal K})\tototi_A M$ on $C_0\subseteq C$, which is
finite over $A_0$. To prove the condition \cite{artin} 4.1.(ii) we can
simply apply the cohomological form of Zariski's theorem on
formal functions, which can easily be extended to
bounded complexes with coherent cohomology.

For condition \cite{artin} 4.1.(iii) we can localize $A_0$ and assume that
the cohomology sheaves of the bounded complex $(\pi_*{\cal
K})\tototi_A M$ are locally free on $\spec A_0$, because
$A_0$ is reduced. Then the complex $(\pi_*{\cal K})\tototi_A
M$ is quasi-isomorphic to a complex with zero differentials,
and the statement is clear.

Condition \cite{artin} 4.1.(i) is straightforward.

Thus we have proven that $\PP$ is a Deligne-Mumford stack.

\section{The weak valuative criterion}
In this section we prove that $\PP$ satisfies the weak valuative criterion for
properness. 
Let $R$ be a discrete valuation ring, $S=\spec(R)$. Let $\eta\in S$ be the
generic point, $s\in S$ the special point. For a finite extension of
discrete valuation rings
$R\subset R_1$, we denote by $S_1, \eta_1, s_1$ the corresponding schemes.

\begin{prp} Let $\X_\eta\to \C_\eta\to \eta$ be a balanced stable fibered
surface. Then there is a finite extension of discrete valuation rings
$R\subset R_1$ and an extension
$$
\begin{array}{ccc}
\X_\eta\times_SS_1 &\subset & \X_1 \\
\down		 &	 & \down \\
\C_\eta\times_SS_1 &\subset & \C_1 \\
\down		 &	 & \down \\
\{\eta_1\} & \subset & S_1, \end{array}
$$
such that $\X_1\to \C_1\to S_1$ is a balanced stable family of fibered
surfaces. The extension is unique up to a unique isomorphism, and its
formation commutes with further finite extensions of discrete valuation
rings.
\end{prp}

\proof We will proceed in steps.

{\sc Step 1: extension of $C$.} Let $X_\eta \to C_\eta\to \eta$ be the
coarse fibered surface, and let $f_\eta:C_\eta\to \m$ be the coarse moduli
morphism. By \cite{fp} and \cite{bm}, there is a complete Deligne - Mumford
stack $\K$. By the weak valuative criterion for $\K$, it follows that there
is a finite extension of discrete valuation rings $R\subset R_1$ and a
unique extension 
\begin{equation}\label{extension-coarse-stable-map}
\begin{array}{ccccc}
C_\eta\times_SS_1 & \subset& C_1 & \stackrel{f_1}{\to} & \m \\
\down		 &	 & \down & &\\
\{\eta_1\} & \subset & S_1 & &\end{array}
\end{equation}
such that $f_1:C_1\to \m$ is a family of Kontsevich stable maps. The
extension \ref{extension-coarse-stable-map} is unique up to a unique
isomorphism and commutes with further base changes.

We now replace $R$ by $R_1$. Let $\pi\in R$ be a uniformizer.

{\sc Step 2: extension of $\C_\eta\to \M $ over $C\sm$.} We may assume that
over
$(C_\eta)\sm$ we have a chart with $V_\eta = U_\eta = (C_\eta)\sm$ and
trivial $\Gamma$: such a chart is compatible with any other chart, so we
may add it to the atlas. We may extend this via $V=U=C\sm$. However the map
$\C_\eta\to \M $ does not necessarily extend over $C\sm$, so we may need a
further base change.

Let $C_i$ be the components of the special fiber $C_s$, and let $\xi_i\in
C_i$ be
the generic points. Consider the localization $C_{\xi_i}$. This is the
spectrum of a discrete valuation ring. By the weak valuative criterion for
$\M$, there is a finite cover $\tilde{C}_i\to C_{\xi_i}$ and a map
$\tilde{C}_i\to \M$ lifting the map on $(C_\eta)\sm$.

We proceed to simplify these schemes $\tilde{C}_i$. Denote $S_n = \spec
R[\pi^{1/n}]$. By Abhyankar's lemma
[SGA I, exp. XIII section 5] we may assume that $\tilde{C}_i$ is \'etale
over $C_{\xi_i}\times_SS_{n_i}$. By descent we already have a lifting
$C_{\xi_i}\times_SS_{n_i} \to \M$. Taking $n$ divisible by all the $n_i$,
we have an extension $C_{\xi_i}\times_SS_{n} \to \M$.

We now replace $R$ by $R[\pi^{1/n}]$. Thus there is an extension $C_{\xi_i}
\to \M$. Note that this extension is unique, since $\M$ is separated. For
the same reason it commutes with further base changes.

By \cite{dj-o}, there is a maximal open set $U\subset C\sm$ with an extension
of the morphism,
$U \to \M$. Since $U$ contains both $(C_\eta)\sm$ and $\xi_i$, we have
$U=C\sm\setminus P$ for a finite set of closed points $P$. By the purity
lemma (Proposition \ref{purity-lemma})
we have that $U = C\sm$, and there is an extension $C\sm \to \M$.

The uniqueness of the lifting in the purity lemma guarantees that this
extension is unique and commutes with further base changes.

{\sc Step 3: extension of $\C_\eta\to \M$ over non-generic nodes.} Let
$p\in \sing(C)$ be a node, and assume $p$ is not in the closure of
$\sing(C_\eta)$.

To build up a chart near $p$, we first choose the \'etale neighborhood
$U\to C$ of $p$, as follows. Let $U_0$ be a Zariski neighborhood such that
$U_0\cap \sing(C) = \{p\}$. We already have a morphism $U_0\to \M$. We take
an \'etale neighborhood
$U$ of $p$ over $U_0$ so that $p$ is a split node on $U$. Thus we can find
elements $s_1,s_2$ in the maximal ideal of $p$, such that $\mm_p=(s_1,
s_2, \pi)$, satisfying the equation 
$s_1s_2=\pi^r$ for some $r>1$. In other words, $U$ is also an \'etale
neighborhood of the closed point $\{s_1=s_2=\pi=0\}$ in $U' = \spec
R[s_1,s_2]/(s_1s_2-\pi^r)$.

We now find a nodal curve $V_0\to S$, which is Galois over $U$, with an
equivariant extension $V_0 \to \M$. Let $V_0'\to U'$ be defined as follows:
$V_0'= \spec R[t_1,t_2]/(t_1t_2-\pi)$, and $s_1 = t_1^r, s_2 = t_2^r$.
There is an obvious {\em balanced} action of the group of $r$-th roots of unity
$\mmu_r$ on $V_0'$ via $ (t_1,t_2) \mapsto (\zeta t_1, \zeta^{-1} t_2)$, and
$V_0'\to U'$ is the 
associated quotient morphism. Let $V_0 = V_0'\times_{U'} U$. The action of
$\mmu_r$ clearly lifts to $V_0$. We write $q_0:V_0 \to U$ for the quotient
map.

Notice that $V_0'$, and therefore $V_0$, is nonsingular. By the purity
lemma we have a lifting $V_0\to \M$ of the morphism $f\circ q_0 : V_0 \to
\m$. Notice that for $g\in \mmu_r$ we have $f\circ q_0 = f\circ q_0 \circ
g$. Thus by the uniqueness of the lifting, the map $V_0 \to \M$ commutes
with the action of $\mmu_r$.

Let $Y_0\to V_0$ be the associated family of stable pointed curves and let
$F_0\subset Y_0$ be the fiber over $p$. The group $\mmu_r$ acts on $Y_0$,
stabilizing $F_0$. Let $\mmu_{r_0} \subset \mmu_r$ be the subgroup acting
trivially on $F_0$, and $\Gamma= \mmu_r/\mmu_{r_0}$ the quotient. Denote by
$Y \to V$ the quotient of $Y_0\to V_0$ by the action of $\mmu_{r_0}$. By
lemma \ref{QuotientIsNodal} we have: $Y \to V$ is a family of stable nodal
curves, $\Gamma$ 
acts on $Y\to V$, and this action is essential.

Let $X_U= Y/\Gamma$.
By \'etale descent, there is an algebraic space $X\to C$ such that $X_U
\simeq X\times_CU$. The existence of an ample sheaf (lemma \ref{ample-sheaf})
guarantees
that $X$ is a scheme.

Thus we have a diagram
$$
\begin{array}{ccccc}
Y	 &\to& X\times_CU & \to & X \\
\down & & \down 	 &	&\down \\
V	&\to& U	 & \to & C \\
\down & & \down &	&\down \\
S	& = &	S	& = &S
\end{array}
$$
giving a chart for $\X \to \C \to S$ near $p$. This is easily seen to
be unique and to commute with further base changes.

{\sc Step 4: extension of $\C_\eta\to \M$ over generic nodes.} Let
$p_\eta\in \sing (C_\eta)$, and let $p\in C_s$ be in the closure of $p_\eta$.
The construction of a chart here is similar to step 2. We choose a Zariski
neighborhood $U_0$ of $p$ such that $U_0\cap \sing(C) = \{p_\eta,p\}$. As
before, we may choose an \'etale neighborhood $U$ over $U_0$ such that
$p_\eta$ is a split node, thus $U$ is \'etale over $\{s_1s_2=0\}$. We may
assume that $V_\eta$ is given by $t_1t_2=0,$ where $t_i^r = s_i$ for an
appropriate integer $r$, with the group $\Gamma = \mmu_r$. Otherwise we can
add a compatible chart with such a $V$ to our atlas. There is an obvious
extension $V_\eta \subset V$ via the same equation $\{s_1s_2=0\}$. We
already have a
lifting $V\setminus \{p\} \to \M$. Since $V$ is normal crossings, the purity
lemma applies, and guarantees that there is a unique equivariant extension
to $V$. It is easy to check that this gives an extension of the chart, which is
essential (as the $\isom$ scheme of stable curves is finite unramified)
and balanced. 

The uniqueness and base change properties are straightforward. \endproof

\section{Boundedness of the category of fibered surfaces}

 The following lemma gives the basic fact behind boundedness:

\begin{lem}\label{finite-fibers} Assume that $S$ is the spectrum of an
algebraically closed 
field. A morphism $C\to \m$ can come from at most finitely many fibered
surfaces $\X\to \C$, up to isomorphism.
\end{lem}

\proof First of all observe that given a family of stable pointed curves
$X'\to C\sm$ there are at most finitely many fibered surfaces $\X\to \C$
whose restriction to $C\sm$ is isomorphic to $X'$; so it is enough to prove
that any given map $f\colon T\to \m$, where $T$ is a smooth curve, comes
from at most finitely many families on $T$.

Over $\bfc$, it is easy to deduce this finiteness result analytically: given a
base point  $t\in T$ and a fixed lifting $\tilde{f}\colon T\to \M$, it is easy
to associate, in a one to one manner, to  any lifting  $\tilde{f'}\colon
T\to \M$ a monodromy representation $\phi_{f'}: \pi_1(T) \to \aut X_t$. Since
$\pi_1(T)$ is finitely generated there are only finitely many such
representations. 

We now give a short algebraic argument using the language of stacks.

 Call $\cal T$ the normalization
of the pullback of $\M$ to $T$; if there is at least one map $T\to \M$
giving rise to $f$ then this induces a section $s\colon T\to {\cal T}$.
It is proved in \cite{vistoli} that this section is \'etale; the
argument is similar to that of the purity lemma: by passing to the strict
henselization one can lift $s$ to the universal deformation space. The image of
this lifting is then evidently \'etale, and therefore $s$ is \'etale.
This strongly restricts the structure of $T$. 

 Call $R$ the fiber
product $T\times_{\cal T} T$; then the \'etale finite groupoid $R\double T$
is a presentation of $\cal T$. Each lifting $T\to \M$ of $f$ comes
from a section $T\to {\cal T}$, so it is enough to show that up to
isomorphisms there are only finitely many such sections. But such a section
$\sigma\colon T\to {\cal T}$ is determined by the pullback $U$of
$T\mapright s {\cal T}$ to $T$ via $\sigma$, and a morphism of groupoids
from $U\times_T U\double U$ to $R\double T$ over $T$. But up to isomorphism
there are only finitely many \'etale coverings $U\to T$ of bounded degree,
because the 
fundamental group of $T$ is topologically finitely generated. For each
$U\to T$ there are only finitely many liftings $U\times_T U\to R$.\endproof

Continuing the proof of Thorem \ref{main-theorem}, we need to show that $\PP$
is proper, and that it admits a projective coarse moduli space $\pp$.

 Since we have poven the weak valuative criterion, in order to prove that $\PP$
is proper, it suffices to show  that it is of finite type. 

The projectivity of $\pp$ then follows, since we have a map $\pp \to \k$ with
finite 
fibers (Lemma \ref{finite-fibers}), which by properness is finite, and
since $\k$ is projective we have 
that $\pp$ is projective.

We begin by adressing the geometric points of $\PP$.
\begin{prp}
There exists a scheme $T$ of finite type over $\bfq$, and a morphism $T \to
\PP$ such that every geometric 
point $\bar{x}\to \PP$ is the image of a geometric point $\bar{x}\to T$ on $T$.
\end{prp}

{\bf Sketch of proof.}
{\sc Step 1: a family of maps to $\m$.} We will build on the fact that $\K$ is
a stack of finite type over $\bfq$. 

Let $T$ be a scheme of finite type over $\bfq$ and $T \to \K$ a morphism
which is surjective on geometric points. We have an associated family
of stable maps:
$$\begin{array}{ccc} C & \to & \m \\
 \dar & & \\
 T. & & \end{array} $$

{\sc Step 2: Straightenning up $C$.}
  Passing to a stratification of $T$, we may assume that 
\begin{itemize}
\item $T$ is nonsingular, and
\item the family $C \to T$ is of locally constant topological type.
\end{itemize}

Let $\tilde{C} \to C$ be the normalization, and let $\Sigma\subset \tilde{C}$
be the locus lying over the nodes $Sing(C/T)$. Since the topological type is
locally constant, $\tilde{C} \to T$ is smooth and $\Sigma \to T$ is
\'etale. Passing to a finite cover of $T$, we may assume that
\begin{itemize}
\item $\tilde{C} = \cup\tilde{C}_i$, where $\tilde{C}_i\to T$ has irreducible
geometric fibers, and 
\item  $\Sigma$ is the union of the images of (disjoint) sections $T \to
\tilde{C}$. 
\end{itemize}

{\sc Step 3: A finite cover with a map to $\M$.}
There exists a projective variety $M$ over $\bfq$ and a finite surjective
morphism $M \to \M$. Let $D = \tilde{C}\times_\m M$. Replacing $T$ by a quasi
finite scheme mapping onto it,  we may assume that the normalization
$\tilde{D}$ of $D$ satisfies:
\begin{itemize}
\item $\tilde{D}\to S$ has geometrically connected smooth fibers;
\item $\tilde{D}\to\tilde{C}$ is flat, branched along the union $\Delta$ of
disjoint sections; and
\item each of the branching sections is either disjoint from $\Sigma$ or
contained in $\Sigma$. 
\end{itemize}

Let $\tilde{C}_0 = \tilde{C} \setminus (\Sigma \cup \Delta)$. Denote by
$\tilde{D}_0$ the inverse image of $\tilde{C}$ in $\tilde{D}$.

{\sc Interlude: discussion of descent data on the open set $\tilde{C}_0$.}
Let $X \to \tilde{D}$ be the family of pointed curves associated with the
morphism $\tilde{D} \to \M$. Let $\bar{t}\to T$ be a geometric point. Fix any
irreducible component $C'\subset  (\tilde{C}_0)_{\bar{t}}$ and an irreducible
component $D'$ above it. Denote the projection map $\gamma:D' \to C'$. Suppose
we have a family of  stable pointed curves $Y\to C'$ such that the associated
moduli map $C' \to \m$ is the given one. Then there is a finite \'etale cover
$\alpha:D'' \to D'$ and an isomorphism $\alpha^*\gamma^*Y' \to
\gamma^*X_{D'}$. We can choose $\alpha$ so that its degree $d$ is
bounded in terms of $\gamma$ and $\nu$ (any bound on the degree of
$\isom_{D'}(X_{D'},X_{D'})$ 
will do). Conversely, any such family of curves $Y \to C'$ comes from descent
data on $\gamma^*X_{D'}$ over $C'$.  Thus it will be useful to put together
such finite \'etale covers $\alpha:D'' \to D'$ . 

Note that any such $D'$ has only finitely many \'etale covers of bounded
degree, since its algebraic fundamental group is topologically finitely
generated.  

{\sc Step 4: Choosing families for basing descent data.}
Replacing $T$ by a surjective quasi finite cover, we may assume
that we have a finite collection of finite maps $D^{(i)} \to \tilde{D}$, which
are \'etale over $\tilde{D}_0$, and such that any connected finite \'etale
cover $\alpha:D'' \to D'$ of  a component of a geometric fiber of
$\tilde{D}_0$, of degree  bounded by $d$ as above, appears in one of
the $D^{(i)}$. We may 
assume as usual that $D^{(i)}\to T$ are all smooth. Let $E^{(i)}$ be the
normalization of $D^{(i)} \times_{\tilde{C}} D^{(i)}$. It admits two
projections 
$p_1, p_2:E^{(i)} \to D^{(i)} \to \tilde{D}$.

{\sc Step 5: Constructing descent data for the open set.}
Descent data on $D^{(i)}_0$ for a family of pointed curves over $\tilde{C_0}$
consist of  a section of the scheme
 $${\mathop{\isom}\limits_{E^{(i)}_0}}(p_1^*X, p_2^*X)$$
 over the open set 
${E^{(i)}_0}$ lying over $\tilde{C}_0$, satisfying the cocycle
condition. Note that by our smoothness assumptions we have that any 
such section extends to a section on $E^{(i)}$. Since $E^{(i)} \to T$ as well
as the $\isom$ scheme are projective over $T$, we have a projective scheme $T'\to T$
(which is in fact finite over $T$) parametrizing such sections. The cocycle
condition puts a closed condition on $T'$. We replace $T'$ by this closed
subscheme. 

Thus, there exists a finite scheme $T'\to T$ parametrizing descent data for
families of curves over geometric fibers of $\tilde{C}_0$. We have the pulled
back family of curves $\tilde{C}'\to T'$, and on the open set $\tilde{C}'_0$ we
have a lifting $\tilde{C}'_0\to \M$. Moreover, all liftings coming from
geometric fibered surfaces appear in this family. 

{}From here on, we can forget about all the $D$'s and $E$'s.

We may again pass to a stratification and assume that $T'$ is nonsingular.

{\sc Step 6: Extending the family over a stack.}
Applying semistable reduction and the purity lemma, it follows that there is a
diagram as follows, unique up to \'etale equivalence: 
$$\begin{array}{ccccc}
& & \tilde{\C}&  \to & \M \\
&\nearrow & \downarrow & & \dar \\
\tilde{C}_0 & \subset & \tilde{C} & \to & \m \end{array} $$
where $\tilde{\C}$ is a Deligne - Mumford stack, smooth over $T$;
$\tilde{\C}\to \tilde{C}$ is the coarse moduli space, which is isomorphic over
$\tilde{C}_0$;  and the map 
$\tilde{\C}\to \M$ is representable, and extends the given lifting on
$\tilde{C}_0$.

There is a subset $T''\subset T'$, which is open and closed, over which the
extension above has the property that  $\tilde{\C}\to \tilde{C}$ is an
isomorphism over $\pi:\tilde{C}\setminus \Sigma$. Similarly, there is a subset
$T'''\subset T''$ over which the stack  $\tilde{\C}$ has a chance to be glued
together over the nodes, namely: for any pair of sections $s_1, s_2: T' \to
\Sigma$ lying over a family of nodes in $C$, the automorphism groups of the
points in $\C$ over them are of constant orders $m_1, m_2$, and $T'''$ is the
locus where $m_1=m_2$ for all such pairs.

{\sc Step 7: Gluing over nodes.}
We might as well replace $T$ by $T'''$. By stratifying $T$, we might as well
assume that  the automorphism groups of the families of curves over $\Sigma$
are locally constant over $T$. We may also assume that the sections of $\Sigma$
lift to a universal deformation space for $\tilde{\C}$, which, we may assume, 
is locally near every node $s$, a cyclic branched cover of $\tilde{C}$ with
Galois group $\Gamma_s$.   

Gluing data for $\tilde{\C}$ over the nodes now consists of a choice of
isomorphism of the objects $s_1, s_2 \in \M(T)$, plus a choice of isomorphism
$\Gamma_{s_1} 
\to \Gamma_{s_2}$. This data is clearly finite over $T$. If we want, we can pass to
an open and closed subset over which the isomorphism makes the gluing stack
$\C$ balanced. The resulting base is $S$. \qed

To finish the proof of the theorem: let $p\in S$ be a point. The family 
$\X_S \to \C_S \to S$ gives an object $\X_p \to \C_p \to p$ in $\PP(p)$.   We
have constructed a 
universal deformation space $V \to \PP$ for any such object, possessing a
family $\X_V \to \C_V \to V$.  Using 
$$\mathop{\isom}\limits_{V \times S} (V\times \X_S  \to
V \times\C_S, \X_V\times S \to \C_V\times S )$$
 we see that there is a nonempty
open set $U$ in $S$ such that all geomertric fibers of $\X_U \to \C_U \to U$
appear in $\X_V \to \C_V \to V$. Continuing by noetherian induction we have a
scheme $V'$ of finite type over $\bfq$ and an \'etale surjective map $V'\to
\PP$, which is what we need.
\qed

\section{Fibered surfaces and Alexeev stable maps}
\subsection{Semi-log-canonical surfaces} 
Semi-log-canonical surfaces are introduced in \cite{ksb} and further studied in
\cite{alex} and  \cite{alexMaps}. We now review their definition. First, let us
define log-canonical pairs.

\begin{dfn} Let $(X,D)$ be a pair consisting of a normal variety $X$ and an
effective, reduced Weil divisor $D$. Denote by $X\sm$ the nonsingular locus of
$X$. We say that $(X,D)$ is a {\em log
canonical pair} if the following conditions hold:
\begin{enumerate}
\item  $(X, D)$ is log-$\bfq$-Gorenstein, namely: for some positive integer
$m$, we assume that the invertible sheaf 
$(\omega_{X\sm}(D))^m$ extends to an invertible sheaf
$(\omega_{X\sm}(D))^{[m]}$ on $X$.  
\item Let $r:Y \to X$ be a desingularization, such that the proper transform of
D 
together with the exceptional locus of $Y\to X$ form a normal crossings
divisor. Call this divisor, taken with reduced structure, $D'$. Then, we assume
that the natural pullback of rational differentials gives a morphism of
sheaves: $r^*(\omega_{X\sm}(D))^{[m]} \to (\omega_{Y}(D'))^m$. That is, a
logarithmic differential on $X$ pulls back to a logarithmic differential on
$Y$. 
\end{enumerate}
\end{dfn}

A complete description of log-canonical singularities of surface pairs is given
in \cite{alexClass}.

Semi-log-canonical surfaces are a natural generalization of the above to the
case of non-normal surfaces. 

\begin{dfn}
Let $X$ be a surface and $D$ a reduced, pure codimension 1 subscheme such that
no component of $D$ lies in $Sing(X)$.  We say that $(X,D)$ is a
semi-log-canonical if the following conditions hold.
\begin{enumerate}
\item $X$ is Cohen Macaulay and normal crossings in codimension 1. Let $X_{nc}$
be the  locus where $X$ is either nonsingular or normal crossings.
\item $(X, D)$ is log-$\bfq$-Gorenstein, namely: for some positive integer $m$,
we assume that the invertible sheaf 
$(\omega_{X_{nc}}(D))^m$ extends to an invertible sheaf on $X$. 
\item Let $X'\to X$ be the normalization, $D'$ the Weil divisor corresponding
of $D$ on $X'$, and $C\subset X'$ the conductor divisor, namely the divisor
where $X'\to X$ is not one-to-one, taken with reduced structure. Then the pair
$(X' , D'+ C)$ is log-canonical.
\end{enumerate}
\end{dfn}

The unique locally free extension of $(\omega_{X_{nc}}(D))^m$ is denoted by
$(\omega_{X}(D))^{[m]}$ 

In case $D$ is empty, we just say that $X$ is semi-log-canonical (or ``has
semi-log-canonical singularities'').

We will use the following lemma about quotients of semi-log-canonical surfaces:
\begin{lem} \label{quot-is-slc}
Let $(Y,D)$ be a semi-log-canonical surface pair. Let $\Gamma\subset \aut(Y,D)$
be a 
finite
subgroup, and let $X=Y/\Gamma$, $D_X = D/\Gamma$.
Then the following
conditions are equivalent:
\begin{enumerate}
\item The quotient $(X,D_X)$ is semi-log-canonical; 
\item
The pair $(X,D_X)$ is log-$\bfq$-Gorenstein.
\end{enumerate}
\end{lem}

The proof is straightforward, see \cite{kollar}.

The reader is referred to \cite{ksb} for the refined notions of semismooth,
semi-canonical 
and semi-log-terminal singularities.

\subsection{Nodal families and semi-log-canonical singularities.}

A family of nodal curves over a nodal curve is always semi-log-canonical:

\begin{lem}
Let $V$ be a nodal curve and $f:Y\to V$ a nodal family. Then $Y$ has
Gorenstein semi-log-canonical singularities. If, furthermore, $D\subset Y$
is a section
which does not meet $\sing(f)$, then $(Y,D)$ is a semi-log-canonical pair.
\end{lem}

\proof Since any family of nodal curves is a Gorenstein morphism, we
have that both $V\to \spec k$ and $Y \to V$ are Gorenstein morphisms,
therefore $Y$ is Gorenstein.

Let $y \in Y$ and $p = f(x) \in V$. It is convenient to replace $Y$ and $V$
by their respective formal completions at $y$ and $p$. The situation falls
into one of the following cases:

\begin{enumerate}
\item Both $p\not\in \sing(V)$ and $y\not\in Sing(f)$. Then $y \in Y_{ns},
$ and there is nothing to prove.
\item $p\not\in \sing(V)$ and $y\in Sing(f)$. Choose a regular parameter
$t$ at $p$. By the deformation theory of a node, we have the description
$$Y \simeq \spf k[[u,v,t]]/(uv- h(t)),$$ with $h(0)=0$. If $h(t)\not\equiv
0$ then we may write
$uv = \mu t^k $ for
some unit $\mu\in \co_{V,p}$, in which case $(Y,y)$ is a canonical
singularity. If $h(t)\equiv 0$ then $(Y,y)$ is a normal crossing point,
which is semismooth.
\item $p\in \sing(V)$ and $y\not\in Sing(f)$. Then $(Y,y)$ is a normal
crossings point again.
\item $p\in \sing(V)$ and $y\in Sing(f)$. We have $$V\simeq \spf
k[[t_1,t_2]]/(t_1t_2) \quad \mbox{ and } \quad Y \simeq \spf O_V[[u,v]]/(uv-
h(t_1,t_2)).$$ We can write $h(t_1,t_2)= h_1(t_1) + h_2(t_2)$, with $h_i(0)=0$.
\begin{enumerate}
\item If neither $h_i$ is $0$, write $h_i(t_i) = \mu_i t_i^{k_i}$, with
$\mu_i$
units, then $(Y,y) $ is a degenerate cusp with exceptional locus a cycle of
rational curves with exactly $\max{(k_1-1,1)} + \max{(k_2-1,1)}$
components. \item If, say, only $h_2\equiv 0$, then we have a degenerate
cusp with $\max{(k_1-1,1)} +2$ exceptional components. \item If $h_1\equiv
0\equiv h_2$, then we have a degenerate cusp with 4 components.
\end{enumerate}
(See [K-SB] \S 4, or [SB] for (4a) - (4c).) \end{enumerate}
The statement about the pair $(Y,D)$ follows easily, since $D$ meets
$\sing(Y)$ transversally at normal crossing points: locally in the \'etale
topology it is isomorphic to the pair $$(\spec
k[t_1,t_2,u]/(t_1t_2),\spec k[t_1,t_2,u]/(t_1t_2,u) ),$$ which is
semi-log-canonical.
\qed

For the following lemma assume that the base field is algebraically
closed. Here we are interested in balanced quotients of families of nodal
curves over nodal base. Studying the resulting singularities becomes easy when
we show that the situation can be deformed, which is shown in this lemma.

\begin{lem}
Let $V$ be the formal completion of a nodal curve and $f:Y\to V$ the formal
completion of a nodal family at a closed point. Let $y\in Y$ be the
closed point and assume $f(y) = p$ is a node. Suppose a finite cyclic group
$\Gamma_y\subset \Aut(Y\to V)$ of order $r$ fixes $y$, acts faithfully on
$V$ stabilizing the two branches of $V$ at $p$ with complementary
eigenvalues (in other words,  the action is balanced).
Then there exists a smoothing
$$\begin{array}{ccc}	 Y& \subset & Y' \\ 			\down & &
\down \\
			 V& \subset & V' \\
			\down & & \down \\
			\spec k & \subset & \spf k[[s]]
\end{array}$$
and a lifting of the action of $\Gamma_y$ to $Y'\to V'$. \end{lem}
{\proof}
\begin{enumerate}
\item Suppose $y\not\in Sing(f)$. We may choose a parameter $u$ along the
fiber so that $Y_y = \spf k[[t_1,t_2,u]]/(t_1t_2)$. It is not hard to
choose $u$ as an eigenvector for $\Gamma_y$. The actions is $$(t_1,t_2,u)
\mapsto (\zeta t_1,\zeta^{-1}t_2,\zeta^a u).$$ This clearly lifts to the
family given by $t_1t_2 = s$.

\item Suppose now $y\in Sing(f)$. We have the local equation $uv=h_1(t_1) +
h_2(t_2)$. It is easy to choose $u,v$ so that $uv$ is an eigenvector.
\begin{enumerate}
\item Suppose neither $h_i$ is $0$. After a change of coordinates we may
assume our local equation is $uv= t_1^{k_1}+ t_2^{k_2}, t_1t_2=0$ (so $r$
divides $k_1+k_2$). We can analyse the action of $\Gamma_y$ via its action
on the fiber $t_1=t_2=0$.
Depending upon whether $\Gamma_y$ stabilizes the branches of $u$ and $v$ or
switches them, the action is either $(u,v) \mapsto (\zeta^au, \zeta^{k_1-a}
v)$ or $(u,v) \mapsto (v, \zeta^{k_1} u)$. In either case this action lifts
to the deformation $uv= t_1^{k_1}+ t_2^{k_2}, t_1t_2=s$. \item Suppose only
$h_2\equiv 0$. We have $uv = t_1^{k_1}, t_1t_2= 0$. Again, the possible
actions are either $(u,v) \mapsto (\zeta^au, \zeta^{k_1-a} v)$ or $(u,v)
\mapsto (v, \zeta^{k_1} u)$.
As before, the action lifts to the deformation. 
\item If we have $uv =
t_1t_2= 0$, the possible actions are either $(u,v) \mapsto (\zeta^au,
\zeta^{k-a} v)$ or $(u,v) \mapsto (v, \zeta^{k} u)$ for some integers $k$
and $a$. For any choice of positive $k_1,k_2$ such that $k_1 - k_2 \equiv k
\mod r$ we have that this action lifts to the family $uv =
t_1^{k_1}t_2^{k_2}; t_1t_2=s$. \end{enumerate}
\end{enumerate}

Now we look at quotients. First the case of a fixed point which lies in teh
smooth locus of a fiber: 

\begin{lem}
Let $V$ be a nodal curve and $f:Y\to V$ a nodal family, let $y\in Y,
y\not\in \sing (f)$. Let $\Gamma\subset \aut Y$ be a balanced finite subgroup,
fixing $y$. Let $q: Y/\Gamma\to X$ be the quotient map and let 
$x=q(y)$. Then $X$ has a semi-log-terminal singularity at $x$. If,
furthermore, $D \subset Y$ is a section through $y$ which is
$\Gamma$-stable, then $(X, D/\Gamma)$ is a semi-log-canonical pair.
\end{lem}

The fact that $X$ has a semi-log-terminal singularity is in [K-SB]
4.23(iii). The statement with $D$ follows as in lemma \ref{quot-is-slc}.

We are left to deal with a fixed point which is a node:

\begin{lem}
Let $V$ be a nodal curve and $f:Y\to V$ a nodal family, let $y\in Y, y\in
\sing (f)$. Let $\Gamma\subset \aut Y$ be a balanced finite subgroup,
fixing $y$. Let $q: Y/\Gamma\to X$ be the quotient map and let $x=q(y)$.
Then $X$ has a Gorenstein semi-log-canonical singularity at $x$. \end{lem}

By lemma ?? it is enough to show that $X$ is Gorenstein, and for this it
suffices to show that the quotient of the smoothing $Y'$ in lemma ?? is
Gorenstein. The quotient variety $Y'/\Gamma$ is clearly Cohen-Macaulay,
therefore it suffices to show
that its canonical divisor class is Cartier. We are in case (2) of lemma
??. In the cases (2a) and (2b) (respectively (2c)), the sheaf $\omega_{Y'}$
is
generated at $y$ by
${{du\wedge dv\wedge dt_2}\over {t_1^{k_1-1}}}$ (respectively, ${{du\wedge
dv\wedge dt_2}\over {t_1^{k_1-1}t_2^{k_2}}}$). 
The generator is
easily seen to be $\Gamma_p$-invariant.
\qed

We thus obtained:

\begin{prp} Let $X\to C$ with sections $s_i:C\to X$ be a coarse fibered
surface, $S_i = Im(s_i)$ and $D=\sum S_i$. Then $(X, D)$ is a
semi-log-canonical pair.
\end{prp}

\subsection{Alexeev stable maps} In \cite{alexMaps}, V. Alexeev defined 
  {\em surface stable maps}, for which he constructed complete moduli
spaces. Our goal here is to compare our moduli of fibered surfaces and coarsely
fibered surfaces with Alexeev's  moduli spaces.

Let $X$ be a reduced, connected projective surface, $D\subset X$ a reduced
subscheme of  codimension 1. Let $M\subset \bfp^r$ be a projective scheme.

A morphism $f:X \to M$ is called a {\em stable map} of the pair $(X,D)$ to $M$,
if
\begin{enumerate} 
\item the pair $(X,D)$ has only semi-log-canonical singularities; in
particular, for some integer $m>0$ the sheaf $(\omega_X(D))^{[m]}$ is
invertible. 
\item For a sufficiently large integer $n$, the sheaf
$(\omega_X(D))^{[m]}\otimes f^*\co_M(mn)$ is ample.
\end{enumerate}

It is easy to see that the property of a morphism $f: X \to M$ being stable is
independent of the choice of the projective embedding of $M$.

Note that, given a stable map $f: X \to M$, one has a well defined
triple of rational numbers: 
$$A =  c_1(\omega_X(D))^2;\quad B= c_1(\omega_X(D))\cdot c_1(f^*\co_M(n));
\quad
C = c_1(f^*\co_M(n))^2.$$

One can define a functor of families of stable maps with fixed invariants $A,B$
and $C$.
This is somewhat subtle, since the sheaf $\omega_X(D)$ is not
invertible, and saturation does not commute with base change. This can be
resolved either by resticting to ``allowable'' deformations on which the
saturation $(\omega_X(D))^{[m]}$ does commute with base change (this is
discussed in an unpublished work by Koll\'ar), or by endowing the surface $X$
with the structure of a 
Deligne-Mumford stack using the log-Gorenstein covers (this has not been
carried out in the litterature). Once the functor is defined, one can look for
a moduli space. The following result of Alexeev gives the answer:

\begin{th}{\bf [Alexeev].} Given rationals $A,B$ and $C$, there is a
Deligne-Mumford stack ${\cal A}l_{A,B,C}(M)$ admitting a projective
coarse 
moduli space ${\bold Al}_{A,B,C}(M)$ for surface stable maps $f:X \to M$ with 
invariants $A,B,C$. 
\end{th}

Now let $\X \to \C \to S$ be a balanced fibered surface, $X \to C \to S$ the
associated 
coarse fibered surface. We have a morphism $C \to \m$, which we can compose
with $X \to C$ and obtain a morphism $f:X \to \m$. In addition to that, we have
$\nu$ sections  $C \to X$. The union of the images of these sectios gives rise
to a divisor $D \subset X$. We have already seen that $(X,D)$ is a
semi-log-canonical pair. We now claim:

\begin{prp} The morphism $f: X \to \m$ is a stable map of the pair $(X,D)$ to
$\m$.  
\end{prp}

\proof Fix a projective embedding $\m\in \bfp^r$. Since $C \to \m$ is a stable
map, we have that $\omega_{C/S}\otimes 
f^*\co_\m(n)$ is ample for $n\geq 3$. Moreover, the line bundle ${\cal L} =
\omega_{X/C}^{[m]}$ defined in lemma  is nef and relatively ample for $X \to
C$. Therefore $\omega_{X/S}^{[m]}\otimes f^*\co_\m(mn)$ is an invertible sheaf,
which is relatively ample for $X\to S$. \qed

It is now easy to see that we have a finite morphism 
$$\PP^{balanced} \to {\cal A}l_{A,B,C}(\m)$$
for a suitable $A,B,C$. Here $\PP^{balanced}$ is teh open substack of balanced
fibered surfaces. The image can be viewed as the stack of
coarse fibered surfaces.

Since the proof of Alexeev's theorem is quite involved, especially showing
boundedness of the stable maps, it is worthwhile seeing that in this particular
case we can deduce the existence of the space of coarse fibered surfaces from
our work so far.

Indeed, the existence of $\PP^{balanced}$ implies that there is a scheme $Z$
and a finite 
surjective morphism $Z \to \PP^{balanced}$. Over $Z$ we have a family of
fibered surfaces, 
and in particular a family of coarse fiberes surfaces $X \to C \to Z$. We have
constructed a relatively ample line bundle $L$ on $X$; we may replace $L$ by a
suitable power and thus assume that $L$ is relatively very ample. After passing
to the frame bundle $Z'$ of sections of $L$, we may assume that the pushforward
of $L$ to $Z$ is a free sheaf $V$. Since $X \to Z$ is flat and embedded in
$\bfp(V)$, we have a morphism $Z \to Hilb$ to a suitable Hilbert scheme. The
projective linear group acts with finite stabilizers, therefore the
quotient is a Deligne Mumford stack, parametrizing coarse fibered surfaces.

It is interesting to look for properties of the morphism $$\PP^{balanced} \to
{\cal 
A}l_{A,B,C}(\m)$$. For instance, it is clearly birational on the closure of the
locus of normal fibered surfaces. As it turns out, there is little more that
one can say: 
this morphism fails in general to be one to one, and moreover, it may be
ramified, even on 
the closure of the locus of normal fibered surfaces. To show this, one simply
needs to produce examples of (deformable) coarse fibered surfaces admitting
incompatible atlases. The following examples are local, but can be easily
globalized.

\begin{example} \label{example1} Here we give an example of two incompatible
balanced charts 
for a coarse  fibered surface. Fix an algebraically closed
field $k$ of characteristic 0, and set $S = \spec k$. Take a smooth
projective curve $W$ over $k$ with an automorphism $s$ of order 2 with a
fixed point $p\in W(k)$, and let $Y_0$ be the curve obtained by attaching
two copies of $W$ at $p$. Consider the action of a cyclic group $\Gamma$ of
order 2 on $X_0$ where a generator acts like $s$ on each copy of $W$; also
let $f$ be the equivariant automorphism of $Y_0$ which acts like $s$ on one
copy of $W$ and as the identity on the other. The point is that $f$ commutes
with $s$, and the automorphism of $Y_0/\Gamma$ induced by $f$ is the
identity.

Let $L$ be  a universal deformation space of $X_0$; $\Gamma$ acts on $L$. Let
$V'$ be a small \'etale neighborhood of 0 in $\bfa^1$, and let 
$\Gamma$ act on $V'$, such that a generator sends $t$ to $-t$. Choose a non
constant $\Gamma$-equivariant map $V' \to L$; this yields a family
$Y'$ of stable curves on $V'$ whose fiber over 0 is exactly $X_0$; $\Gamma$
also acts on $Y$. Let $V$ be the union of two copies of
$V'$ glued at 0; there are two families $Y_1$ and $Y_2$ on $V$ obtained by
attaching two copies of $Y'$ at $X_0$, one using the identity, the other
using $f$. These will not be isomorphic in general, not after an
automorphism of $V$, and not even after going to an \'etale neighborhood of 0
in $V$. 

Set $U = V/\Gamma$ and call $X$ the union of two  copies of
$Y'/\Gamma$ along $X_0/\Gamma$; we claim that the two generically fibered
surfaces
$X_1 = Y_1/\Gamma \to V/\Gamma = U$ and $X_2 = Y_2/\Gamma \to V/\Gamma = U$
are canonically isomorphic to $X\to U$. In fact the structure sheaf
$O_{X'_i}$ fits into an exact sequence
$$0\longrightarrow O_{X'_i} \longrightarrow O_{Y'}\times
O_{Y'}\longrightarrow O_{X_0}\longrightarrow 0,$$ where the map
$O_{Y'}\times O_{Y'}\longrightarrow O_{X_0}$ is the difference of the two
projections for $i = 1$, and the difference of one projection with the
other projection twisted by $f$ for $i = 2$. Now take invariants; we get an
exact sequence $$0\longrightarrow
O_{X'_i/\Gamma}\longrightarrow O_{Y'/\Gamma}\times O_{Y'/\Gamma}\longrightarrow
O_{X_0/\Gamma}\longrightarrow 0;$$
where the two maps
$O_{Y_1/\Gamma}\times O_{Y_2/\Gamma}\to O_{X_0/\Gamma}$ are equal, because
$f$ induces the identity on $O_{X_0/\Gamma}$.

So $(U,Y_1\to V,\Gamma)$ and $( U,Y_2\to V,\Gamma)$ are incompatible charts
for the surface $X\to U$. If we call $\overline U$ the union of two
disjoint copies of $U$, $\overline V$ the union of two disjoint copies of
$V$ and $\overline Y\to V$ the family of stable curves which coincides with
$Y_1$ on one copy and with $Y_2$ on the other one, we even obtain an
example $( \overline U,\overline Y\to V,\Gamma)$ of a chart that is not
compatible with itself. This is not surprising, because we are using the
\'etale topology, and so $U$ is not in general embedded in $C$; it could
not happen if we were working over {\bf C} with the analytic topology.
\end{example}

\begin{example} \label{example2} We now give an example of two incompatible
charts 
on a fibered surface over $S=\spec k[\epsilon]/(\epsilon^2)$, which
coincide modulo $\epsilon$. We then globalize this example to a
complete curve. This implies, in particular, that the map $\PP
\to Al(\m)$ is not always unramified.

Let $F = \spec k[u,w] / (uw)$. We will use $F$ as a constant fiber in
charts $Y \to V$, so that $Y = F \times V$. Specifically, consider the
following two curves:  $V_0 = \spec k[z,t,\epsilon] / (zt, \epsilon^2)$
and $V_\epsilon = \spec k[z,t,\epsilon] / (zt - \epsilon, \epsilon^2)$.
Define $Y_0 = F \times V_0$ and $Y_\epsilon = F \times V_\epsilon$. We
have obvious morphisms $Y_0 \to V_0 \to S$ and $Y_\epsilon \to
V_\epsilon \to S$. The fibers over $\spec k \subset B$ are clearly
isomorphic, but $Y_0$ and $Y_\epsilon$ are clearly non-isomorphic.

We now define an action of $C_6$, the cyclic group of order 6, on
these schemes. Let $\zeta$ be a primitive sixth root of 1. We choose a
generator of $C_6$ and call it $\zeta$ as well; define its action as follows:

$$(u,w,z,t,\epsilon) \mapsto (\zeta^3u, \zeta^2w, \zeta z,
\zeta^{-1}t, \epsilon)$$

This clearly defines an action of $C_6$ on both $Y_0 \to V_0 \to S$
and $Y_\epsilon \to V_\epsilon \to S$.

Let $C_0 = V_0 / C_6$. Explicitly, $$ C_0 = \spec k [z^6, t^6,
\epsilon]/(z^6t^6, \epsilon^2).$$ Similarly let $C_\epsilon
=V_\epsilon / C_6$. It is easy to see that we similarly have
$$ C_\epsilon = \spec k [z^6, t^6,
\epsilon]/(z^6t^6, \epsilon^2).$$ Thus we can use the notation $$C = V_0 /
C_6 = V_\epsilon / C_6.$$

We denote $X_0 = Y_0 / C_6$ and $X_\epsilon = Y_\epsilon / C_6$.
Clearly $X_0\to C \to S$ and $X_\epsilon\to C \to S$ are coarse fibered
surfaces over $S$. Our main claim is that they are isomorphic.

We will produce an isomorphism of coarse fibered surfaces by chosing
isomorphisms over suitable open sets, showing that these glue together
over a large open set, and arguing that an isomorphism on the large
open set must extend. Let us first work out the open sets.

We denote by $U_0$ the localization of $Y_0$ at $u$, \ldots , $T_\epsilon$
the localization of $Y_\epsilon$ at $t$. Explicitly,

\begin{eqnarray*}
 U_0 & = & \spec k [ u, u^{-1}, z,t,\epsilon]/(zt, \epsilon^2) \\
 W_0 & = & \spec k [ w, w^{-1}, z,t,\epsilon]/(zt, \epsilon^2) \\
 Z_0 & = & \spec k [ u, w, z,z^{-1},\epsilon]/(uw, \epsilon^2) \\
 T_0 & = & \spec k [ u, w, t,t^{-1},\epsilon]/(uw, \epsilon^2) \\
 U_\epsilon & = & \spec k[u,u^{-1},z,t,\epsilon]/(zt-\epsilon,\epsilon^2) \\
 W_\epsilon & = & \spec k[w,w^{-1},z,t,\epsilon]/(zt-\epsilon,\epsilon^2) \\
 Z_\epsilon & = & \spec k [ u, w, z,z^{-1},\epsilon]/(uw, \epsilon^2) \\
 T_\epsilon & = & \spec k [ u, w, t,t^{-1},\epsilon]/(uw, \epsilon^2)
\end{eqnarray*}

These open sets are clearly stable under the action of $C_6$.

The expressions above immediately give isomorphisms $Z_0\simeq
Z_\epsilon$ and  $T_0\simeq
T_\epsilon$. These canonically induce isomorphisms  $Z_0/C_6\simeq
Z_\epsilon/C_6$ and  $Z_0/C_6\simeq
Z_\epsilon/C_6$.

Let us address the open sets $U_0$ and $U_\epsilon$.  These are
clearly nonisomorphic, but we claim that their quotients by the
subgroup of two elemets are isomorphic. Indeed, note that $\zeta^3$
 acts trivially on $u, u^{-1}$. Therefore
\begin{eqnarray*}
U_\epsilon / C_2 &\simeq& \spec k[u, u^{-1}] \times V_\epsilon/C_2 \\
 & = & \spec k[u, u^{-1}] \times \spec
k[z^2,t^2,zt,\epsilon]/(zt-\epsilon, \epsilon^2) \\
&=& \spec k[u, u^{-1}]\times \spec k[z^2, t^2,\epsilon]/(\epsilon^2)\\
&\simeq& \spec k[u, u^{-1}] \times V_0 / C_2 \\
&\simeq& U_0.
\end{eqnarray*}

Taking the quotient by the residual $C_3$ action, this isomorphism
clearly induces an isomorphism $$U_\epsilon / C_6 \simeq U_0 / C_6.$$

Note that on the intersections $U_\epsilon/C_2\cap Z_\epsilon/C_2$ and
$U_0/C_2 \cap Z_0 / C_2$ the isomorphism above coincides with the
restriction of the isomorphism $Z_\epsilon/C_2\simeq Z_0/C_2$. Indeed, both
are given by identifiying the variables $u,z^2$ and $\epsilon$.
Therefore these isomorphisms coincide on $U_\epsilon/C_6\cap
Z_\epsilon/C_6$ as well. A similar situation is obtained by
replacing $Z_\epsilon$ by $T_\epsilon$.

We now address the open sets $W_\epsilon$ and $W_0$ in a similar
manner, noting that this time $C_3$ acts trivially on the variable
$w$. Thus
\begin{eqnarray*}
W_\epsilon / C_3 &\simeq& \spec k[w, w^{-1}] \times V_\epsilon/C_3 \\
 & = & \spec k[w, w^{-1}] \times \spec
k[z^3,t^3,zt,\epsilon]/(zt-\epsilon, \epsilon^2) \\
&=& \spec k[w, w^{-1}]\times \spec k[z^3, t^3,\epsilon]/(\epsilon^2)\\
&\simeq& \spec k[w, w^{-1}] \times V_0 / C_3 \\
&\simeq& W_0.
\end{eqnarray*}

Again, it is easy to check that this isomorphism agrees on the
interstections with the open sets $Z_\epsilon$ and $T_\epsilon$.
We note that all these isomorphisms are isomorphisms over the curve
$C$, which coincide with the identity modulo $\epsilon$.

Let $X_\epsilon'$ be the open subscheme which is the complement of the
image of the point $u=w=z=t=0$ on $Y_\epsilon$; define $X_0'$ in an
analogous way. Thus far we have obtained  an isomorphism
$\phi':X_\epsilon' \to X_0'$. The closure $\Sigma\subset
X_\epsilon\times_CX_0$ of the graph of $\phi'$ is a subscheme
supported along the  graph of the identity $(X_\epsilon)_{\spec k} =
(X_0)_{\spec k}$. Therefore the projections $p_\epsilon:\Sigma \to
X_\epsilon $ and $p_0:\Sigma \to X_0$ are finite. Note also $p_\epsilon$
and $p_0$ are isomorphic along a dense open set. Since the schemes
$X_0$ and $X_\epsilon$ satisfy Serre's condition $S_2$, these
projection admit sections, giving rise to a morphism $\phi:X_\epsilon
\to X_0$, which by the same reason is an isomorphism.
\end{example}

\end{document}